\documentclass[hidelinks,onefignum,onetabnum]{siamart220329}
\usepackage{geometry}
\usepackage{lipsum}
\usepackage{amsfonts}
\usepackage{graphicx}
\usepackage{epstopdf}
\usepackage{algorithmic}
\ifpdf
  \DeclareGraphicsExtensions{.eps,.pdf,.png,.jpg}
\else
  \DeclareGraphicsExtensions{.eps}
\fi

\newsiamremark{remark}{Remark}
\newsiamremark{hypothesis}{Hypothesis}
\crefname{hypothesis}{Hypothesis}{Hypotheses}
\newsiamthm{claim}{Claim}

\headers{Adaptive Partition of Unity Interpolation Method with Moving Patches}{A. Heryudono and M. Raessi}

\title{Adaptive Partition of Unity Interpolation Method with Moving Patches\thanks{This is a preprint version.
\funding{This work was funded by the National Science Foundations DMS-2012011.}}}

\author{Alfa Heryudono\thanks{Department of Mathematics, University of Massachusetts Dartmouth, Dartmouth, MA, 02747, USA
  (\email{aheryudono@umassd.edu}).}
\and Mehdi Raessi\thanks{Department of Mechanical Engineering, University of Massachusetts Dartmouth, Dartmouth, MA, 02747, USA 
  (\email{mraessi@umassd.edu}).}
  }
\usepackage{amsopn}

\ifpdf
\hypersetup{
  pdftitle={Adaptive Partition of Unity Interpolation Method with Moving Patches},
  pdfauthor={A. Heryudono and M. Raessi}
}
\fi

\newcommand{\ulx}{\ensuremath{\underline{x}}}

\DeclareMathOperator{\Tr}{Tr}

\begin{document}

\maketitle

% REQUIRED
\begin{abstract}
The adaptive partition of unity interpolation method, introduced by Aiton and Driscoll, using Chebyshev local interpolants, is explored for interpolating functions with sharp gradients representing two-medium problems. For functions that evolve under vector fields, the partition of unity patches (covers) can be shifted and resized to follow the changing dynamics of local profiles. The method is tested for selected 1D and 2D two-medium problems with linear divergence-free vector fields. In those cases, the volume fraction in each patch contributing to volume conservation throughout the domain can be kept in high accuracy down to machine precisions. Applications that could benefit from the method include volume tracking and multiphase flow modeling.
\end{abstract}

% REQUIRED
\begin{keywords}
partition of unity method, free boundary problems, volume-preserving technique
\end{keywords}

% REQUIRED
\begin{MSCcodes}
65D99, 65M70, 76M25 
\end{MSCcodes}

\section{Introduction}
The need for highly accurate interpolants for functions where the region of high activity is shifting or moving in time appears in many practical applications such as shape deformations and multiphase flows. Depending on the methods used, such functions are usually sampled on structured/unstructured grids, meshes, or scattered points. Initially, if the methods permit, the underlying grids are typically concentrated in some regions to capture localized profiles such as peaks and sharp gradients. As the profiles move in time to some other locations in the domain, the grids usually need to be redistributed or resampled accordingly. Hence, the accuracy and quality of the interpolant can be maintained. The interpolant is typically utilized not just for data interpolation but also for finding derivatives and computing integrals. 

\begin{figure}[htbp]
\begin{minipage}{\linewidth}
\begin{center}
\begin{minipage}{0.35\linewidth}
\scalebox{0.225}{\includegraphics{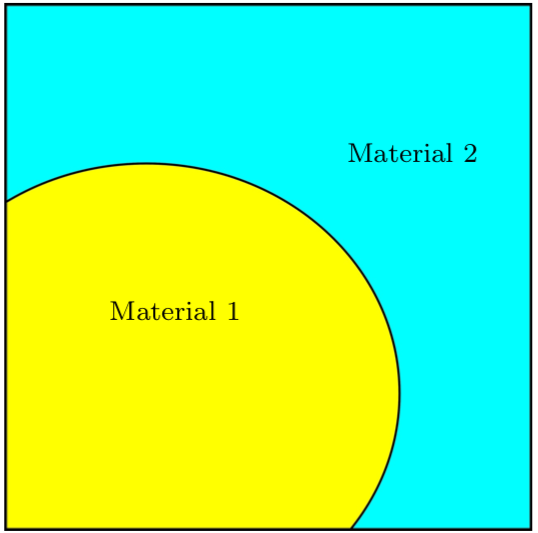}}
\end{minipage}
\hspace{10pt}
\begin{minipage}{0.35\linewidth}
\begin{displaymath}
f(\underline{x})=\left\{\begin{matrix}
1, & \underline{x} \in \mathrm{Material~1}\\ 
0, & \underline{x} \notin \mathrm{Material~1}
\end{matrix}\right.
\end{displaymath}
\end{minipage}
\end{center}
\end{minipage}
\caption{An illustration of two materials/mediums layout along with the interface. $f(\underline{x})$, usually called an indicator function, is commonly modeled with a discontinuous function or a function with sharp gradient.}
\label{fig:twomat}
\end{figure}

In this work, we are interested in cases where a one- or two-dimensional time-dependent function $f(\underline{x},t)$ defined on a rectangular domain $\Omega$ has sharp gradients representing a two-phase medium scenario; an example is shown in \cref{fig:twomat}, where $ f(\underline{x},t) = 1$ (representing material/medium 1), on a particular localized region inside the domain and zero everywhere else. The locations where the sharp gradients occur are called the interface, which separates the two mediums. With a given vector field $\underline{u}(\underline{x})$, with $\underline{x} = (x,y)$ in 2D, for example, one might be interested in studying the dynamics of the free boundary (the shape of the interface) and the conservation of several quantities associated with the mediums/materials. To be specific, with our function described above, $\int{f} d\Omega$ measures the volume of medium $1$. In problems where volume conservation is expected, the value of the integral must stay the same throughout the simulation. Typically, one can rarely use the grid points used at time $t=0$ for the initial condition $f(\underline{x},0)$ since the profile of the function might look different later. Ideally, grids should be distributed following the moving front.  

Computational strategies and numerical methods for dealing with functions representing a two-phase medium evolved under a vector field consist of two challenges. (a) A highly accurate interpolant for $f$ in space that captures the region of high activities is needed, and (b) distributing the grids/points/cells marching together with the time-evolving interface while ensuring the quality of the interpolant does not change. There has been a large body of publications in this field since the early development of numerical methods for partial differential equations. For example, it is at the heart of multiphase flow simulations using the Volume-of-Fluid (VOF) method.

This paper explores ways to mitigate challenges (a) and (b) by using the partition of unity method with non-stationary covers. To approximate functions with sharp gradients in this work, we focus on utilizing an adaptive partition of unity (APU) method based on Chebyshev polynomial interpolants introduced by Aiton, and Driscoll \cite{aiton_adaptive_2018,aiton_adaptive_2019,aiton_preconditioned_2020}. The method is able to adaptively (using bisection techniques in alternating dimensional directions) construct a highly accurate global interpolant down to machine precision with ease in 1D, 2D, and 3D. Their codes are available in \textsc{MATLAB} and use Chebfun \cite{Driscoll2014} as a backend. Both make rapid numerical prototyping with Chebyshev technologies enjoyable with minimal effort. 

During its construction, the interpolant also creates overlapping patches (covers) of the computational domain and its local interpolants. For functions with high activities around a particular region (e.g., around the interface) in the domain, more patches are concentrated there and less anywhere else. To follow the dynamics of the free boundary or the interface, we allow the patches to move along with the flow field as long as they meet certain conditions. For updating local function values in the patches, the resizing and shifting of the covers usually only involve scaling due to changes in covers' sizes. Hence, the accuracy of the global interpolant is not affected.

The paper is organized as follows: In Section 2, we briefly introduce the partition of unity interpolation method. Section 3 describes how patches can resize and shift due to a vector field and how function values are updated. The algorithm is briefly discussed in Section 4. Numerical cases of two-phase problems in 1D and 2D under linear divergence-free vector fields that show conservation of volume down to machine precision are provided in Section 5. A discussion section concludes our work. 

\section{Partition of Unity Method in a Nutshell}
\label{sec:pum}

The basic idea of the partition of unity approach is to break the domain into several pieces (patches), approximate the function in each subdomain (patch) separately, and then blend the local approximations together using smooth, local weights that sum up to one everywhere on the domain. Let $f(\underline{x})$ be defined as a smooth scalar multivariate function defined on a closed domain $\Omega\subset \mathbb{R}^n$, with boundary $\partial \Omega$ where $\underline{x}=(x_1,x_2,...,x_n)$. In a partition of unity method, the global approximation $\tilde{f}(\underline{x})$ to the function $f(\underline{x})$ is constructed as a weighted sum of local approximations $\tilde{f}_j(\underline{x})$ on overlapping patches $\Omega_j$, $j=1,\ldots, N_p$. That is,
\begin{align}
\tilde{f}(\underline{x})=\sum_{j=1}^{N_p} w_j(\underline{x})\tilde{f}_j(\underline{x}).
\label{eq:PUM}
\end{align}
where $w_j$, $j=1,\ldots,N_p$ are weight functions.
The patches $\Omega_j$ need to form a cover of the domain in the sense that
\begin{align}
\bigcup_{j=1}^{N_p}\Omega_j\supseteq \Omega.
\end{align}
The partition of unity weight functions 
$w_j$ are non-negative, compactly supported on $\Omega_j$ and satisfy
\begin{align}
\sum_{j=1}^{N_p}w_j(\ulx)=1,\quad \forall \ulx\in\Omega.
\end{align}

In the particular case where rectangular patches are used, the local approximations $\tilde{f}_j(\underline{x})$ are Chebyshev approximations on the local tensor product grid. In the best efficiency scenario, the local approximant is computed only once and replicated for each patch. However, several templates of local approximants can also be precomputed and stored in a lookup table to be used by any patch that needs them. In order to capture high degrees of localization, patches can be adapted to reflect the profile of the solutions. Hence, the method offers greater flexibility to adjust local approximants by refining/coarsening patches (h-type version) or raising/lowering the degree of local interpolants (p-type version).

Calculus operations such as differentiation and integration of the global approximant \eqref{eq:PUM} can be done term by term (then summed up) in a straightforward way. This is handy when gradients, divergence, and volumes must be computed. The calculus operations should be accurate locally at patch levels and globally throughout the domain.

Constructing a PU interpolant for a function representing two-medium problems usually results in patches clustered in the interface's neighborhood. Those patches containing local interpolants with high-degree Chebyshev polynomials are needed in regions with steep gradients. Otherwise, the interpolant may not capture the function properly there, i.e., under resolution can most likely happen. Indeed, the locations of those patches can also be utilized to determine the approximate shape of the interface. This information could potentially be helpful for the purpose of tracking or reconstructing interfaces, although we have not used the Chebyshev PU interpolation method solely for that purpose here yet. Other regions away from the interface can have patches with larger sizes with lower degrees of local interpolant.

\begin{figure}[htbp]
\begin{minipage}{\linewidth}
\hspace{-25pt}
\scalebox{0.375}{\includegraphics{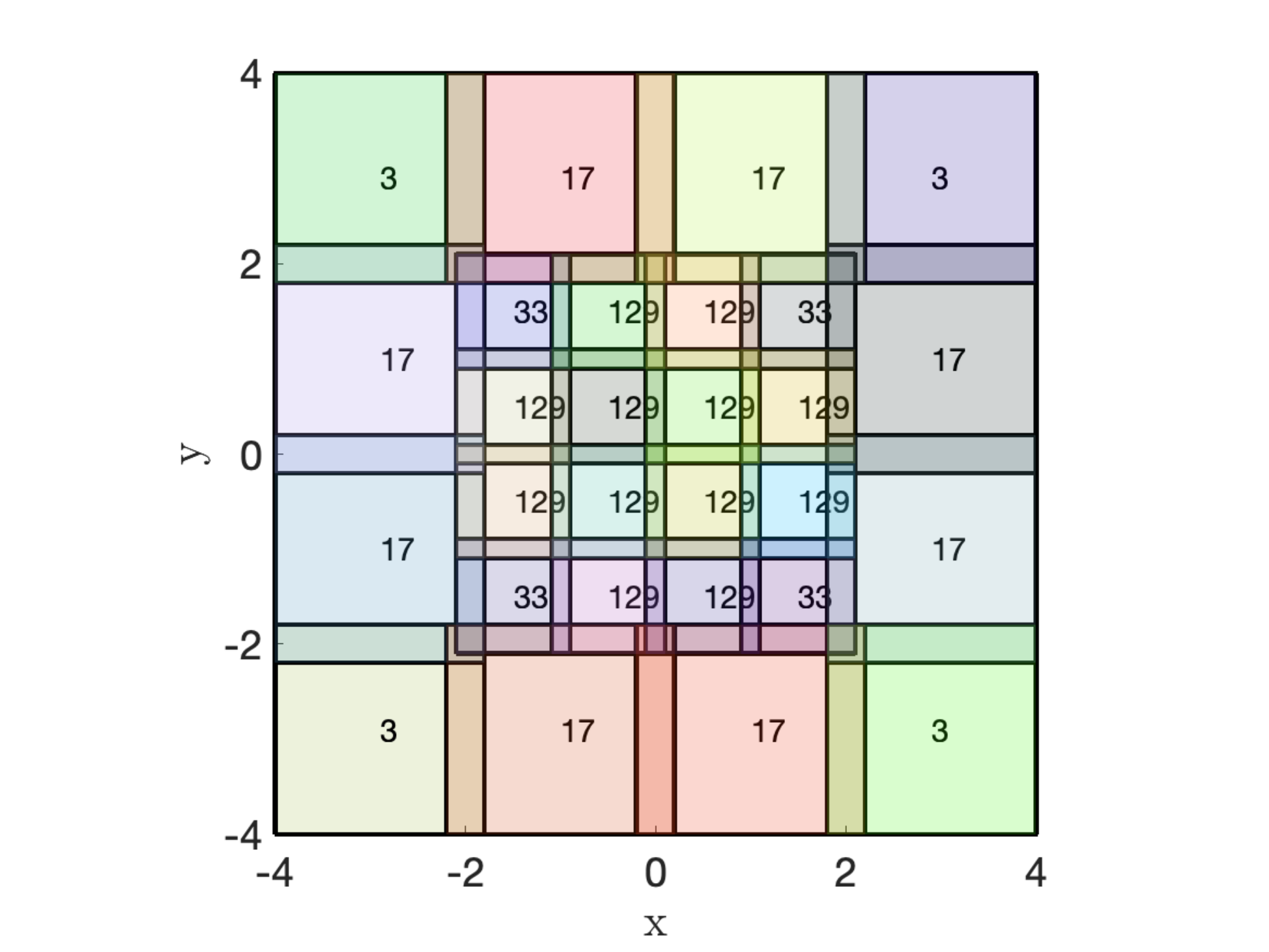}}
\hspace{-35pt}
\scalebox{0.375}{\includegraphics{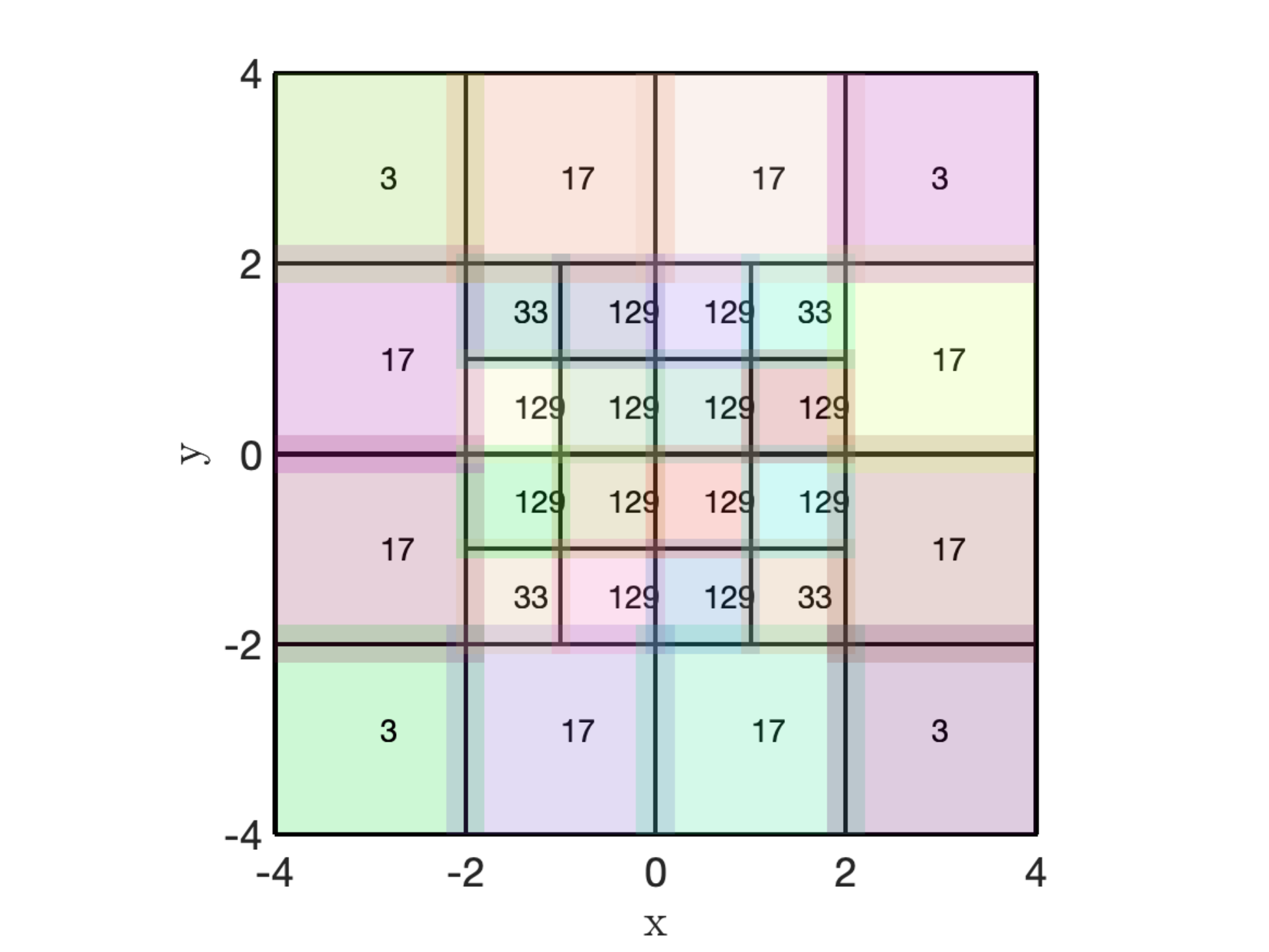}}
\end{minipage}
\caption{The layout of patches of the function $f$, systematically constructed using bisections, along with their overlapped regions. The numbers show the maximum degrees of Chebyshev polynomials used in each patch. The figure on the right shows the zones: non-overlapping boxes that define the patches.}
\label{fig:pillustration}
\end{figure}

As an illustration in two dimensions, \cref{fig:pillustration} shows the layout of the rectangular patches for a function $f(x,y) = \tfrac{1}{2}(1 + \tanh(5(1-(x^2+y^2))))$ in the domain $[-4,4] \times [-4,4]$, with polynomials degrees ranging from $3$ to $129$ with 10\% overlap. For visualization purposes, those rectangles can be colored based on any information that might be useful for users if needed. Aiton and Driscoll use a semitransparent pastel style coloring so that the overlap regions can also be observed. 

\begin{figure}[htbp]
\begin{minipage}{\linewidth}
\scalebox{0.415}{\includegraphics{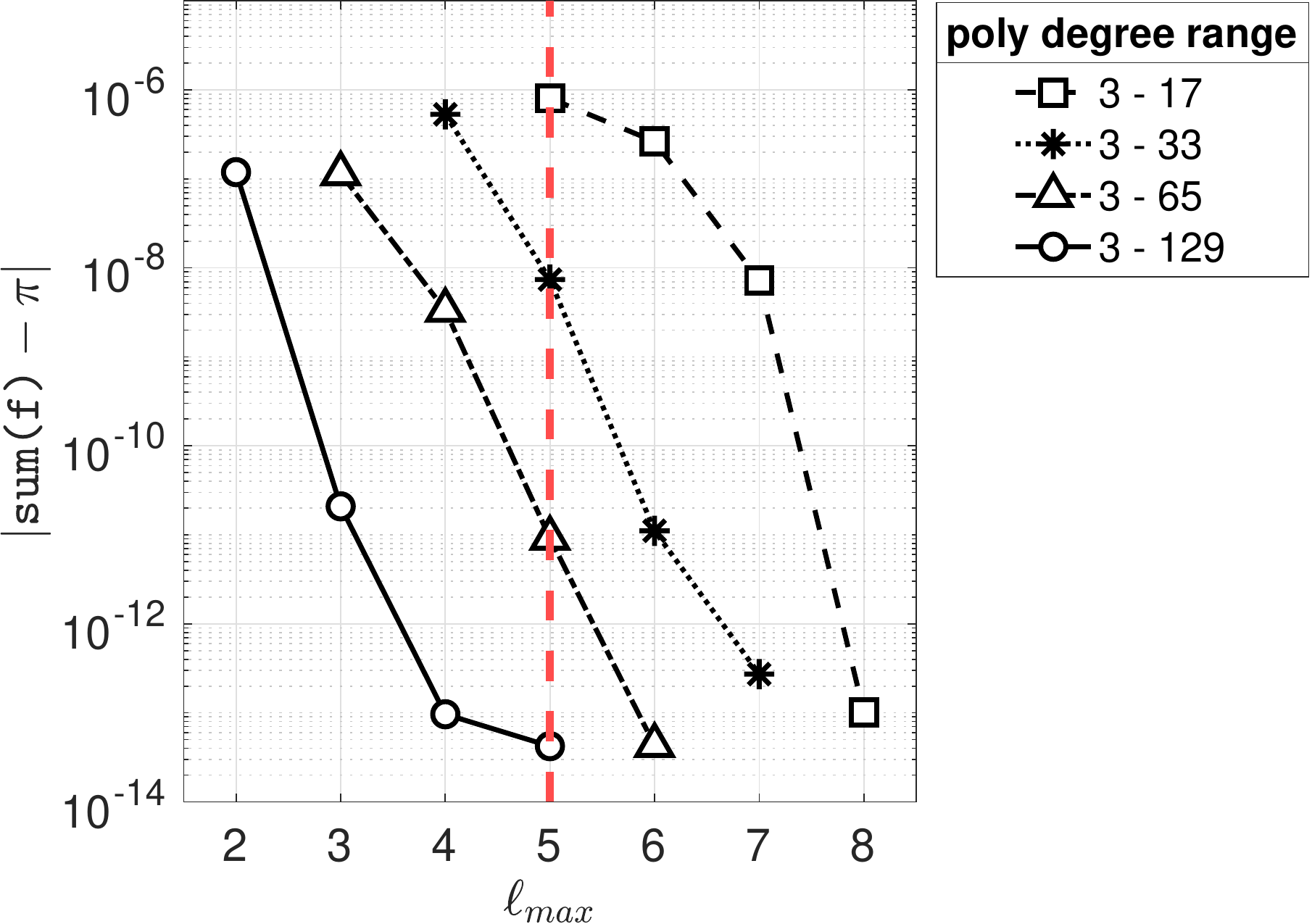}}
\hspace{2pt}
\scalebox{0.395}{\includegraphics{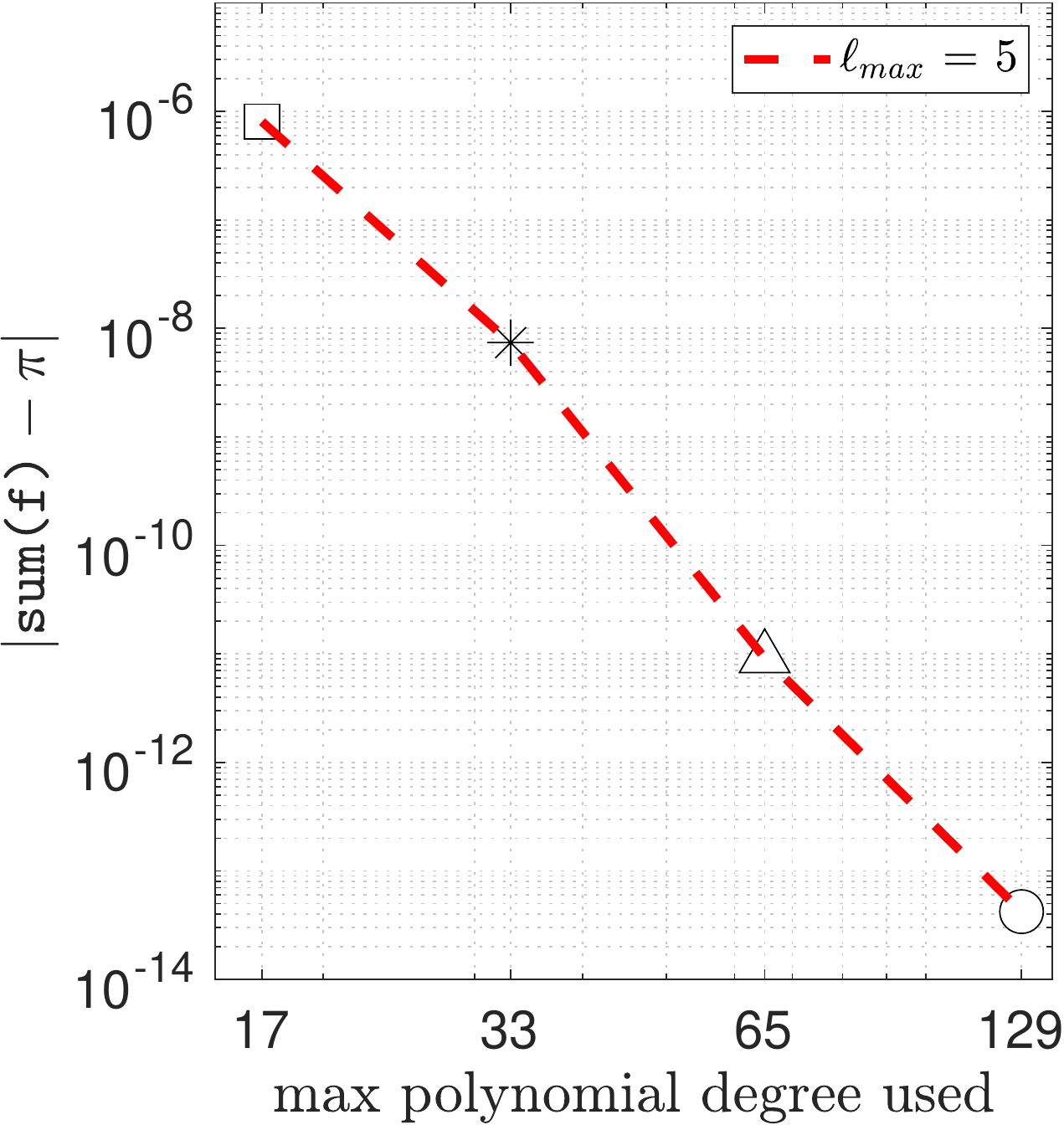}}
\end{minipage}
\caption{The function $f(x,y) = \tfrac{1}{2}(1 + \tanh(100(1-(x^2+y^2))))$ is constructed using an APU interpolant. The \texttt{sum(f)} $\approx \int f d\Omega$ is computed and the value then is compared to $\pi$. Left: error $|\texttt{sum(f)} - \pi |$ convergence trend (\texttt{semilogy} scale) as a function of $\ell_{max}$ for varying polynomial degree range used by the interpolant in each patch. In this case, the maximum subdivision of the patch/zone in either $x$ or $y$ direction can go down to $2^{-\ell_{max}}.  $Right: the convergence (\texttt{loglog} scale) by fixing the $\ell_{max}$ (red dash line on the left figure) and plot it as a function of maximum polynomial degree used in each patch.}
\label{fig:numexconv}
\end{figure}

With their codes, to construct the interpolant for the function above, one may execute the following command in \textsc{MATLAB}:

\begin{verbatim}

f = PUchebfun(@(x,y) 0.5*tanh(5*(1-(x.^2+y.^2)))+0.5,[-4 4;-4 4]);

\end{verbatim}
With the backend of Chebfun, the construction usually finishes in a fraction of a second but could be longer for more complicated functions. Similar commands are available in 1D and 3D. Once the interpolant is constructed, one can use it to interpolate or find operators involving derivatives (e.g., Laplacian, gradients, divergence) and integrals. For example, one can execute \texttt{sum(f)} (i.e., $\int f d \Omega$) to compute the integral of $f$ on the domain. In the case above, the integral of $f$ approximates the area of a unit circle. When one increases the slope of the $\tanh$ from $5$ to say $50$, the integral gives about $14$ digits accuracy in a fraction of a second. One can see the power of adaptive partition of unity representation for constructing a global interpolant with spectral accuracy.

Plotting the error convergence trend of the adaptive method, though not as straightforward as in the fixed grid method, can be done from different scenarios. This is due to the different range of polynomial degrees and the level of subdivisions used during the construction of the interpolant. We provide two example plots. \cref{fig:numexconv} shows the convergence plot of $|\texttt{sum(f)} - \pi|$ using the APU interpolant for the same hyperbolic tangent with the slope raised again to $100$. The left figure of \cref{fig:numexconv} shows the error convergence trend by allowing the degree of Chebyshev polynomials to vary in either $x$ or $y$ direction in each patch/zone but restricting the maximum level of subdivision of the patches/zones down to $2^{-\ell_{max}}$. On the other hand, the rightmost figure of \cref{fig:numexconv} shows when one wants to fix the maximum level (say $\ell_{max} = 5$) and then plot the convergence along that particular vertical line. In the end, the goal of the adaptive method is to systematically use a combination of subdivisions and a range of polynomial degrees to achieve certain error tolerance with minimal/no user intervention. Moreover, the class constructor $f$ reveals the underlying fields and structures useful for numerical experiments or prototyping. Information about local patches, weights, Chebyshev series coefficients, and others is provided as part of the data structure. The papers \cite{aiton_adaptive_2018,aiton_adaptive_2019,aiton_preconditioned_2020} along with the code explain everything users need to know about this method. 

\section{Moving patches}
\label{sec:movpatch}

When the profile of the function changes with respect to a vector field, patches may not stay stationary anymore. Ideally, patches with high-degree polynomials are always closely following the interface. Alternatively, one can always rebuild the patches by initiating the interpolant constructor at every step. However, this process should probably be done not too often or under some specific conditions only to reduce computational time. In simple cases where only translation and pure strain problems are considered, rebuilding patches are usually unnecessary.
Indeed, in translation and pure-strain vector fields, a particular patch is responsible for a specific subdomain of the function at all times, regardless of whether it is stretched or shrunk. Hence ``material transfer'' that requires donor-acceptor techniques between patches can be avoided. 

In the current work, several constraints or conditions exist to consider when allowing rectangular patches to move following a linear divergence-free flow field. These conditions are common in computer graphics for simulating shape deformations of objects \cite{von_funck_vector_2006}. Some of these restrictions might be possibly removed in future work.

\begin{enumerate}
\item The movement of patches is based on the movement of the vertices of the patches.
\item The shape of the patches must stay rectangular when stretched or shrunk. Cases with patch rotations are left for future study. The rectangular condition is needed for this work, which uses a structured grid. However, the shape can be somewhat arbitrary for other partition of unity frameworks with unstructured grids. 
\item Patches cannot detach or merge. Overlapped regions should remain intact too.
\item Their relative positions must always stay the same, i.e., neighbors for life, even when their sizes change or positions move.
\item The number of patches stays the same throughout the simulation. Some on-the-fly refinement/coarsening strategies will be left for future study.
\item Boundaries of the domain can be fixed or freely moved with the flow. Fixed boundaries can cause a patch to stretch with a higher proportion in one direction.  
\end{enumerate}
Although the movement of those patches may look all over the place during the simulation, they are moving in an orchestrated way with respect to those constraints. 

Given a vector field, the time-dependent position vector of the vertices of the patches can be obtained by solving systems of ODEs. However, as a proof of concept, we provide the position vector function or utilize MATLAB and Chebfun high-order built-in ODE solvers to maintain high accuracy in both space and time when needed. This may not be realistic in applications where $\underline{u}$ may come from other flow solvers. At least, with the APU interpolant, the spatial approximation is accurate, and the error of the conserved quantities may be due to the accuracy of the ODE solver used.

For the 2D problem, in each patch, the types of linear divergence-free velocity field $\underline{u}(x,y) = \left < u_1(x,y), u_2(x,y) \right >$ we are considering in this work are classified in \cref{tab:typesvel}.
\begin{table}[htbp]
\footnotesize
\caption{Types of velocity fields. All $c$'s and $\omega$ are constants. $\omega$ is also referred to as angular velocity.}\label{tab:typesvel}
\begin{center}
 \begin{tabular}{|c|c|c|} \hline
  Type & \bf $u_1(x,y)$ & \bf $u_2(x,y)$ \\ 
  \hline
  Translation & $c_1$ & $c_2$ \\
  Pure Strain & $c x + c_1$ & $-c y + c_2$ \\ 
  Angular Deformation & $c y + c_1x$ & $c x-c_1y$ \\ 
  Pure Rotation & $-\omega y + c x$ & $\omega x - c y$ \\ 
  \hline
 \end{tabular}
\end{center}
\end{table}
The divergence-free conditions provide desired properties for shape deformation such as no path line self-intersections (local or global) in 4D space-time domain \cite{theisel_topological_2005} and the deformation is volume-preserving \cite{davis_introduction_1979}.

In terms of system of linear ODEs, the velocity field can be written as
\begin{displaymath}
\begin{bmatrix}
u(x,y)\\
\\
v(x,y)
\end{bmatrix}
=
\begin{bmatrix}
\frac{dx}{dt}\\
\\
\frac{dy}{dt}
\end{bmatrix}
= \begin{bmatrix}
& & \\
& A & \\
& & \\
\end{bmatrix}
\begin{bmatrix}
x \\
\\
y \\
\end{bmatrix}
+
\begin{bmatrix}
 \\
b \\
\\
\end{bmatrix},
\end{displaymath}
with $\Tr(A) = 0$ ensuring the divergence-free condition. Entries of $A$ and $b$ consist of $c, c_1, c_2, \omega$ depending on the type of the vector fields.

\begin{figure}[htbp]
\begin{minipage}{\linewidth}
\begin{minipage}{0.45\linewidth}
\centering
\scalebox{0.65}{\includegraphics{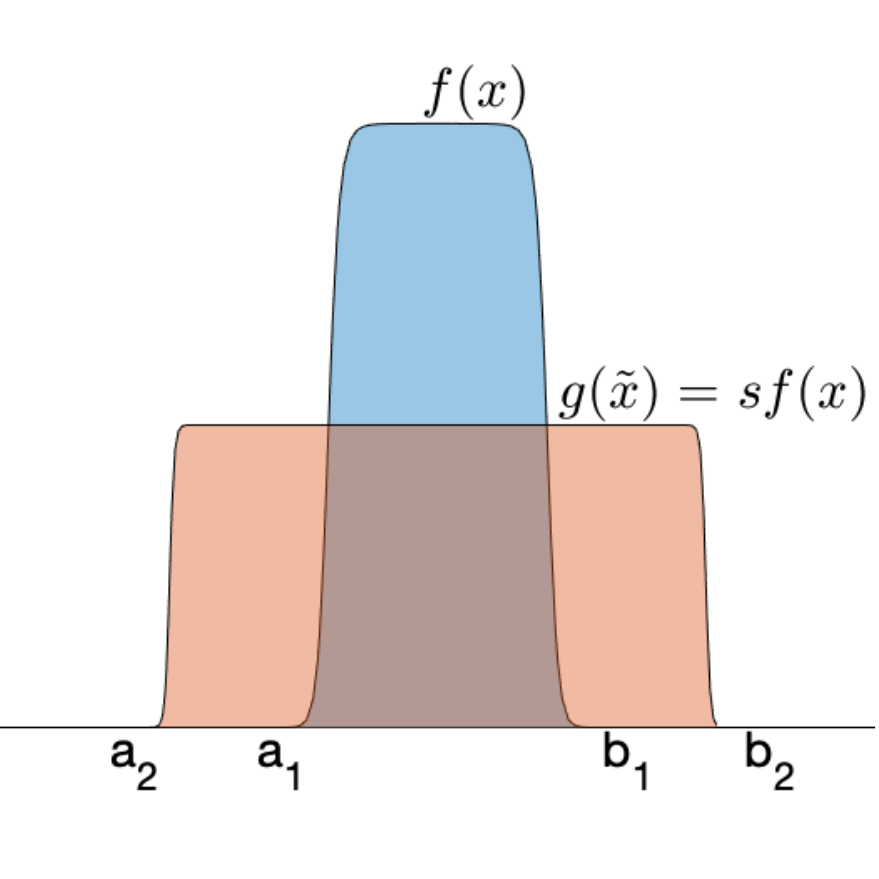}}
\end{minipage}
\begin{minipage}{0.45\linewidth}
\centering
\scalebox{0.55}{\includegraphics{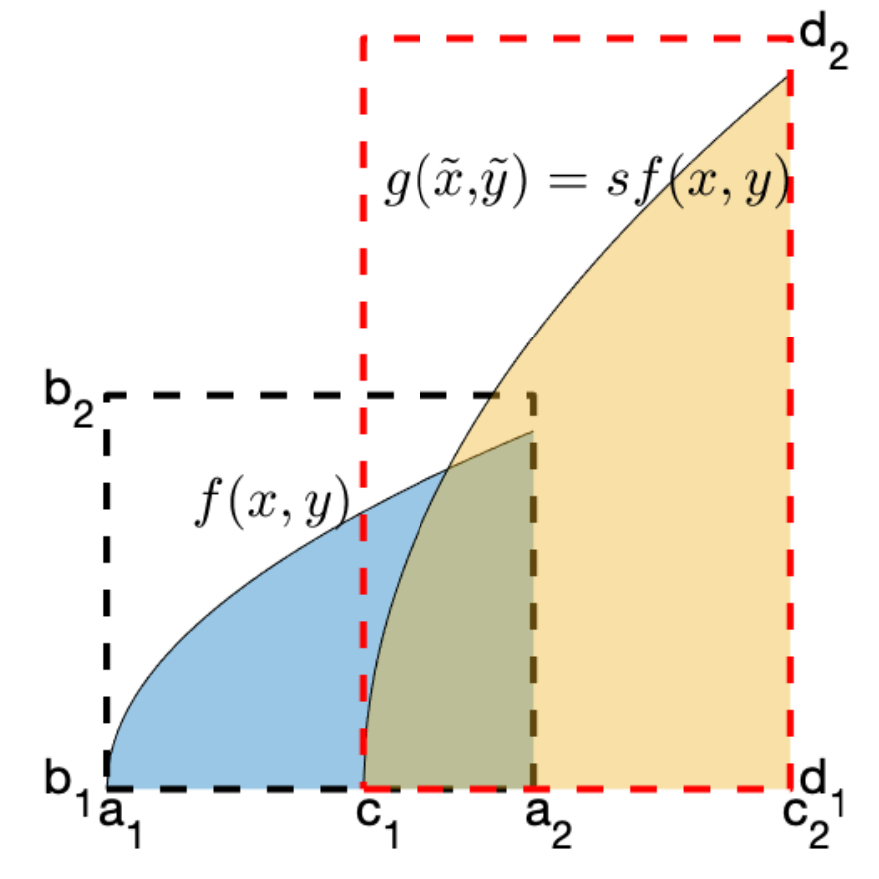}}
\end{minipage}
\end{minipage}
\caption{After a patch is resized, the function $f$ is scaled by a constant $s$ to keep the volume constant. Deforming the function $f$ this way means shape can change but volume and mass stay constant.}
\label{fig:volscaling}
\end{figure}

We can set a condition for ``deforming" the function $f$ on $\Omega$ or $f_j$ or $w_j f_j$ when a patch $\Omega_j$ is shifted or resized. Without loss of generality, we just use the notation $f$ for both global and local functions. Note that the model function we are using here, for example, hyperbolic tangent in 1D, is Lipschitz continuous on $[a,b]$ or on any patch interval. It means that there is a constant $C$ such that $|f(x) - f(y)| \leq C |x - y|$ for all $x,y \in [a,b]$. Locally, the product of the function with the partition of unity weight function in any particular patch is also Lipschitz continuous. \cref{fig:volscaling} shows an illustration in 1D. When a patch of interval $[a_1,b_1]$ is resized to $[a_2,b_2]$, the goal is to scale the function $f(x)$, $x \in [a_1,b_1]$ to a new function $g(\tilde{x})$ such that $g(\tilde{x}) = s f(x)$, $\tilde{x} \in [a_2,b_2]$ with
\begin{displaymath}
\int_{a_2}^{b_2} g(\tilde{x}) d\tilde{x} = \int_{a_1}^{b_1} f(x) dx
\end{displaymath}  
to keep the area/volume unchanged. In this case, $s = \frac{|b_1 - a_1|}{|b_2 -a_2|}$. In 2D, as illustrated on the right figure on \cref{fig:volscaling}, the 
$s = \frac{|a_2 - a_1|}{|c_2 - c_1|} \frac{|b_2 - b_1|}{|d_2 - d_1|}$. Since the change of the variable from $x$ to $\tilde{x}$ is a linear map, $s$ is the product of ratio of the size of the patch before and after resizing it in each dimension. In other words, deforming the function this way means shape can change but volume and mass stay constant.

For approximating local interpolants with the Chebyshev series and their integrals, each patch interval can be scaled to $[-1,1]$. The following two theorems and their proof, see \cite{trefethen_approximation_2013}, can be used.

\begin{theorem}[Chebyshev Series]\label{thm:csthm}
If $f$ is Lipschitz continuous on $[-1,1]$, it has a unique representation as a Chebyshev series,
\begin{displaymath}
  f(x) = \sum_{k=1}^\infty a_k T_k(x),
\end{displaymath}
which is absolutely and uniformly convergent. The coefficients are given for $k \geq 1$ by the formula
\begin{displaymath}
  a_k = \frac{2}{\pi} \int_{-1}^{1} \frac{f(x) T_k(x)}{\sqrt{1-x^2}} dx,
\end{displaymath}
and for $k=1$, by the same formula with factor $2/\pi$ changed to $1/\pi$.
\end{theorem}
The adaptive partition unity method monitors the decay of the coefficients in each dimension to decide to split the domain to construct patches and chop the series at the desired $n$. The highest order of polynomial it can use in each dimension is $129$, though users can modify the APU code if needed. 

Once the Chebyshev series coefficients are available, then they can be used to calculate, for example, an integral in each patch. The theorem below describes the formula.
\begin{theorem}[Integral of a Chebyshev series]\label{thm:icsthm}
The integral of a degree $n$ polynomial expressed as a Chebyshev series 
 \begin{displaymath}
 \int_{-1}^1 \sum_{k=1}^n c_k T_k(x) dx = \sum_{k=0, k\;\textnormal{even}}^n \frac{2c_k}{1-k^2}.
 \end{displaymath}
\end{theorem}
The total integral for the whole domain is the sum of the term-by-term integral of \cref{eq:PUM} contributed from each patch. In each patch, the 2D and 3D implementation of the APU constructor is the extension of the 1D version via tensor product. Note that $a_k$ and $c_k$ coefficients are already precomputed and scaling the function $f$ with $s$ will leave them unmodified.

\section{Algorithm}
\label{sec:alg}

The algorithm of the partition of unity method with moving patches can be described in \cref{alg:pump}. We choose to use a hyperbolic tangent (also commonly used in THINC approach \cite{xiao_simple_2005}) with a steep slope as a template for modeling a two-medium function $f$. 

\begin{algorithm}
\caption{Adaptive partition of unity method with moving patches.}
\label{alg:pump}
\begin{algorithmic}
\STATE{Initiate a function $f$ constructor with adaptive partition of unity method.}
\STATE{Compute $v_0 =$ the total volume of $f$ on $\Omega$.}
\STATE{Set starting time $t=T_\textnormal{init}$, a time step $\Delta t$, and a final time $T_\textnormal{fin}$ for the simulation.}
\WHILE{$t <  T_\textnormal{fin}$}
\FOR{each patch $\Omega_i$}
\STATE{Get all the vertices coordinates $vertcoords(\Omega_i(t))$ at time $t$.}
\STATE{Record the volume $vol(\Omega_i(t))$ at time $t$.}
\STATE{Solve the system of ODEs at the time interval $[t,t+\Delta t]$ based on the vector field information to obtain $vertcoords(\Omega_i(t + \Delta t))$.}
\STATE{Record the new volume $vol(\Omega_i(t + \Delta t))$ at time $t + \Delta t$.}
\STATE{Scale the function values with the scaling factor $s = \frac{vol(\Omega_i(t))}{vol(\Omega_i(t + \Delta t))}$.}
\ENDFOR
\STATE{Update all patch vertices with all new information at $t+ \Delta t$.}
\STATE{Compute $v =$ the total volume of $f$ and measure the error $|v - v_0|$.}
\STATE{update the time $t=t+\Delta t$.}
\ENDWHILE
\end{algorithmic}
\end{algorithm}

Note that if the vector field $\underline{u}(\underline{x})$ is provided beforehand, the trajectory of the vertices can be precomputed. We do that here in most of our numerical experiments to show the method as a proof of concept. However, in practice, the flow field is obtained from the flow solver, and an analytical formula is unavailable. Although we have not done it here, the for-loop step in \cref{alg:pump} can be done in parallel.

\section{Numerical Experiments}
\label{sec:numtest}

All our experiments here can be classified as velocity-field-based shape deformation problems using linear divergence-free vector fields. Our numerical experiments are carried out in \textsc{MATLAB} 2021a on a standard workstation equipped with a quad-core Intel CPU with 16GB RAM. Aiton provided the Github repository where the PUchebfun codes can be downloaded \cite{puchebfungithub}. Chebfun \cite{Driscoll2014} software can be downloaded from its website. Both packages should be added to the \textsc{MATLAB} path. We provide a simple 1D code explaining the idea in Experiment 1. As proof of concepts, in all our numerical experiments, all the flow fields are provided. Hence, the timing of most computations provided here can be done in a few seconds. In the future, we will couple this with a flow solver.

\subsection{Experiment 1 (expanding a 1D blob)}

We begin by evolving the initial condition at time $t=0$ given by the function
\begin{displaymath}
f(x) = \tfrac{1}{2}\tanh(\upsilon(1-x^2))+\tfrac{1}{2},
\end{displaymath}
with slope $\upsilon=100$ on the interval $-4 \leq x \leq 4$ under the flow field $u(x) = cx$, where $c=1/\upsilon$. In terms of two-phase medium terminology, $f$ can be thought as an indicator function where the region $-1 \leq x \leq 1$ represents material $1$ and the rest of the region represents material $2$. The simulation is stopped at $t=1.6$ to avoid the effect of boundaries. 

First, we create the APU constructor for the function by executing as well as computing its volume (the area under the curve) at $t=0$ with the following commands.

\begin{verbatim}

dom = [-4 4];
f = PUchebfun(@(x) 0.5*tanh(100*(1-x.^2))+0.5,dom)
v0 = sum(f)

\end{verbatim}
This creates a global partition of unity interpolant $f$ with 22 patches. The degree of Chebyshev polynomials in those patches varies between $3$ and $129$. As expected, patches with the highest degrees of $129$ are located around $x=\pm1$ where the interfaces exist.

\begin{figure}[htbp]
\begin{minipage}{\linewidth}
\begin{minipage}{0.32\linewidth}
\begin{center}
\scalebox{0.35}{\includegraphics{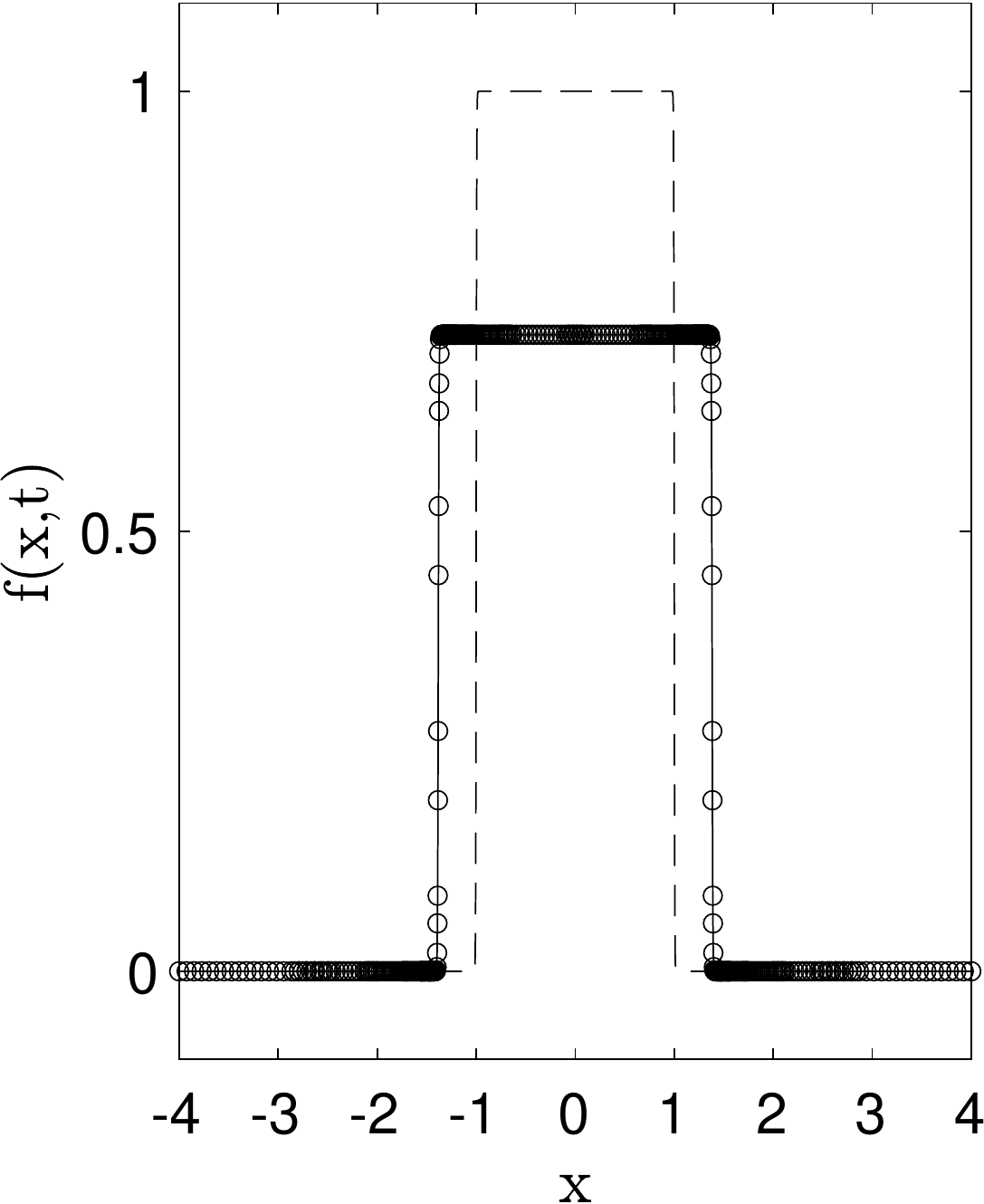}}
\end{center}
\end{minipage}
\begin{minipage}{0.32\linewidth}
\begin{center}
\scalebox{0.35}{\includegraphics{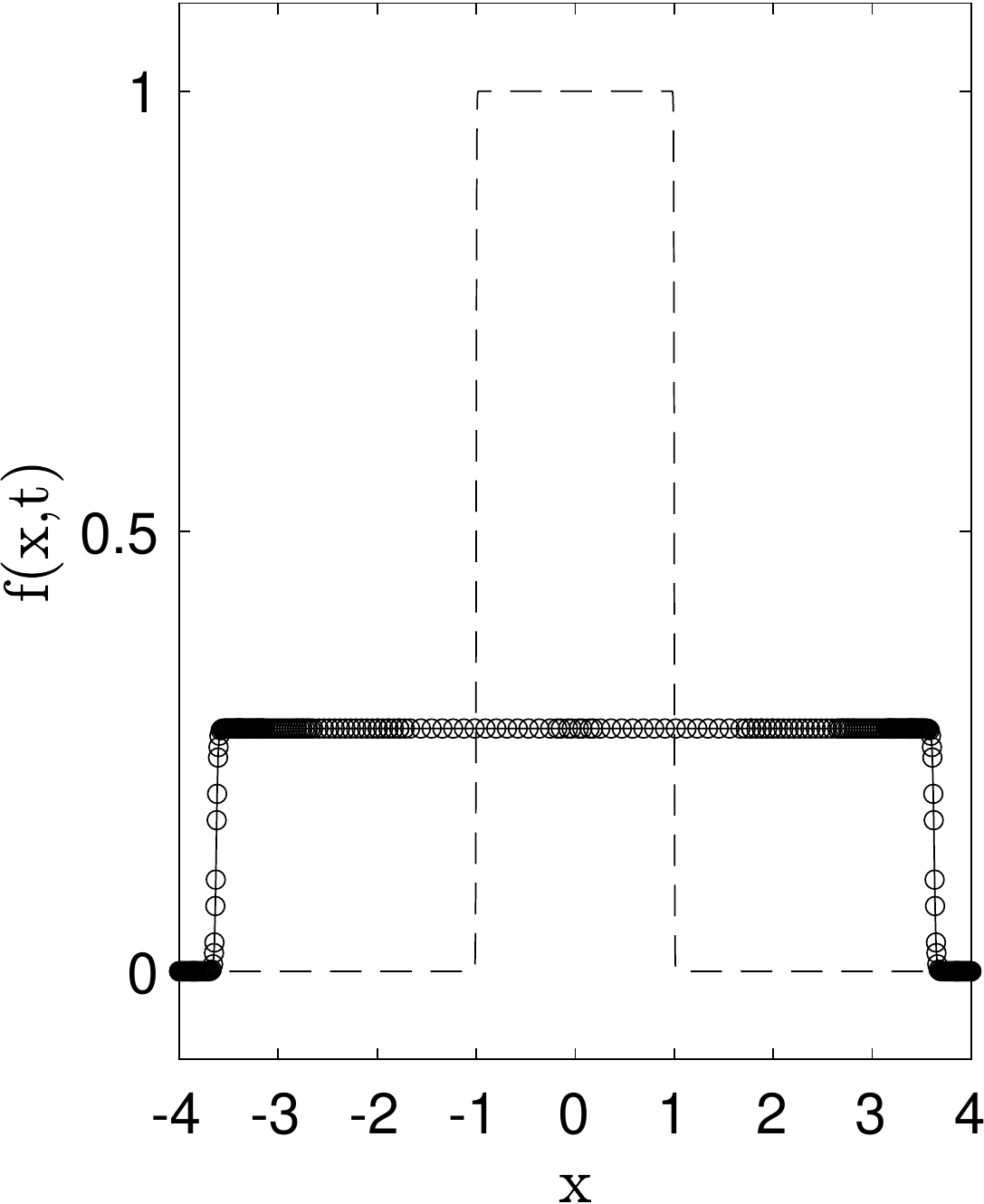}}
\end{center}
\end{minipage}
\begin{minipage}{0.32\linewidth}
\begin{center}
\scalebox{0.35}{\includegraphics{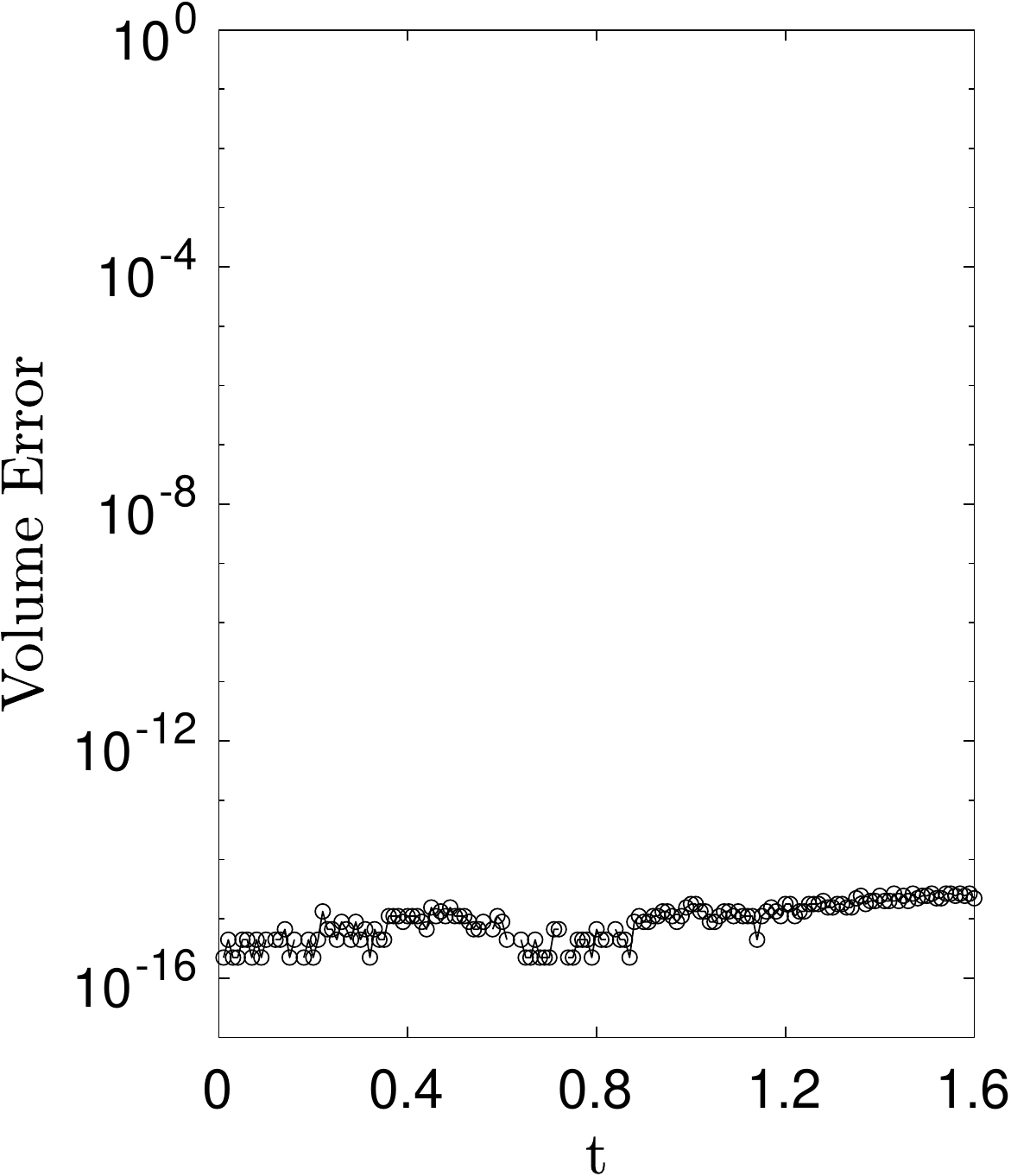}}
\end{center}
\end{minipage}
\end{minipage}
\caption{The evolution of the $\tanh$ from $0 \leq t \leq 1.6$. The dash line is the function at $t=0$. The left most figure is the profile at time $t=0.8$ and the middle figure is at time $t=1.6$. The right most figure shows the absolute error between the volume at time $t$ with the one computed at initial time. The deformation does seem to preserve the sharp gradient with no smoothing or energy minimization needed.}
\label{fig:numex1}
\end{figure}

In terms of time, under the flow field $u(x)=c x$, the endpoints of each patch are evolving as $x(t) = x_0e^{ct}$. The only points that are not moved on purpose are the interval endpoints $x=\pm 4$. As shown in \cref{fig:numex1}, the simulation results in stretching material $1$ with decreasing height to maintain a constant volume. Following the scaling technique provided in Section 3, the decreasing height is due to the multiplication with a ratio of the patch size before and after resizing. The volume at every step can be quickly computed because all coefficients are already precomputed.

Since the analytic trajectory is already known beforehand, we can use it immediately in the code. The rough template of the \textsc{MATLAB} code, with $\Delta t = c = \frac{1}{\upsilon}$ without the plot, is provided as the following:

\begin{verbatim}

Np = length(f.leafArray);
[domlen,newdomlen] = deal(zeros(Np,1));
c = 1/100; x = @(t) exp(c*t); 
Tinit = 0; dt = c; Tfin = 1.6;
for t=Tinit:dt:Tfin
   for i=1:Np
    domlen(i) = diff(f.leafArray{i}.domain);
    flag = f.leafArray{i}.domain > dom(1) & ...
           f.leafArray{i}.domain < dom(end);
    f.leafArray{i}.domain(flag) = f.leafArray{i}.domain(flag)*x(t);
    f.leafArray{i}.zone(flag) = f.leafArray{i}.zone(flag)*x(t);
    newdomlen(i) = diff(f.leafArray{i}.domain);
    ratlen = domlen(i)/newdomlen(i);
    f.leafArray{i}.values = f.leafArray{i}.values*abs(ratlen);
   end
   v = sum(f);
   abs(v-v0)
end

\end{verbatim} 
The template of the code is relatively similar in 2D and 3D since one can access the data structure of the PUchebfun easily. Note that one can instead also pass a function handle obtained from ODE solvers replacing $x(t)$. The volume error, defined as the difference between the initial volume and calculated volume at any time, remains very near machine zero, as shown in the right-most figure in \cref{fig:numex1}. The highest volume error was no greater than $3 \times 10^{-15}$. When one uses ODE solvers in replacing $x(t)$, one might expect the error in the volume plot to be at the order of accuracy of the time-stepping used.

\subsection{Experiment 2 (reversible blob in 1D)}

\begin{figure}[htbp]
\begin{minipage}{\linewidth}
\begin{minipage}{0.32\linewidth}
\begin{center}
\scalebox{0.35}{\includegraphics{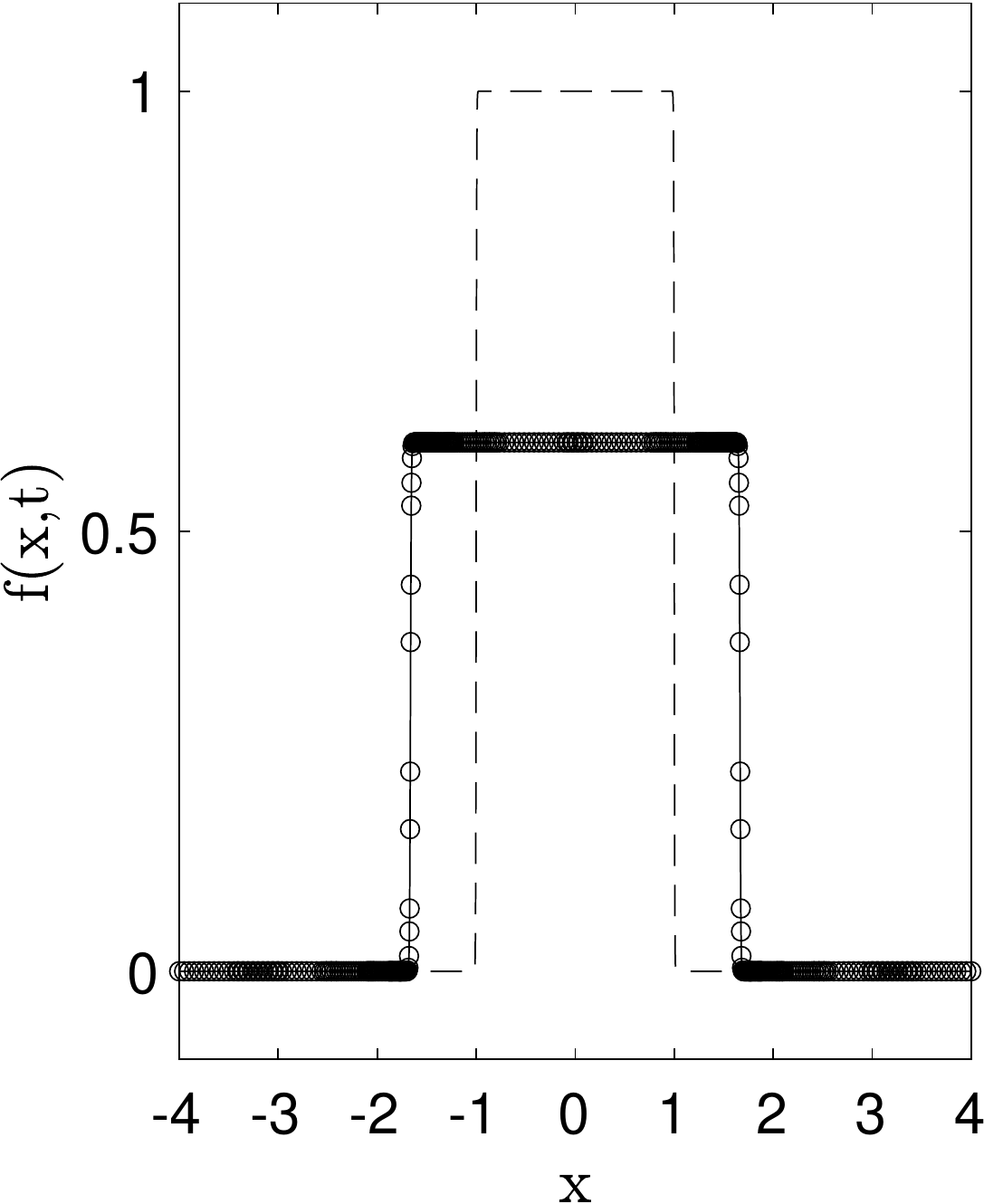}}
\end{center}
\end{minipage}
\begin{minipage}{0.32\linewidth}
\begin{center}
\scalebox{0.35}{\includegraphics{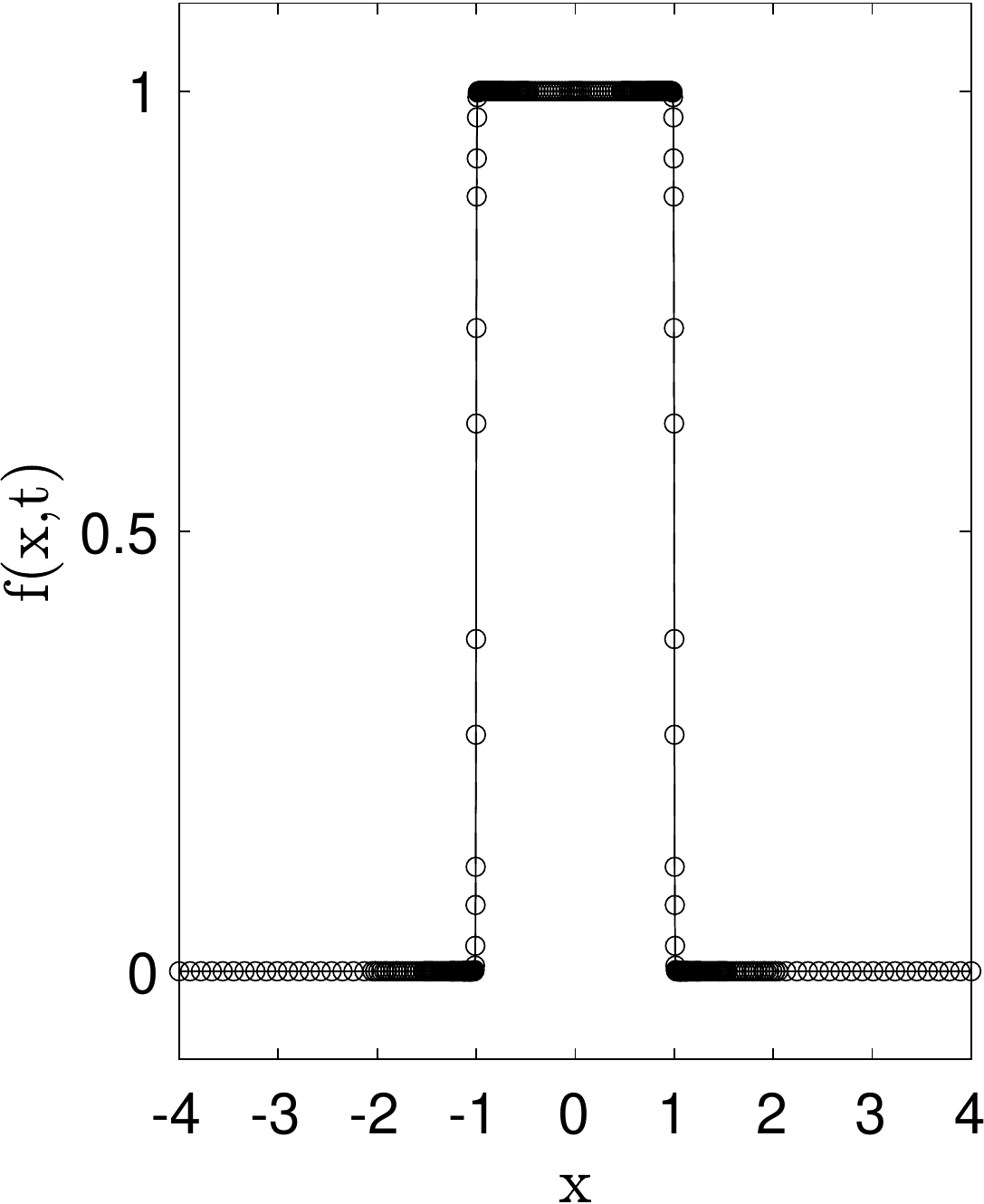}}
\end{center}
\end{minipage}
\begin{minipage}{0.32\linewidth}
\begin{center}
\scalebox{0.35}{\includegraphics{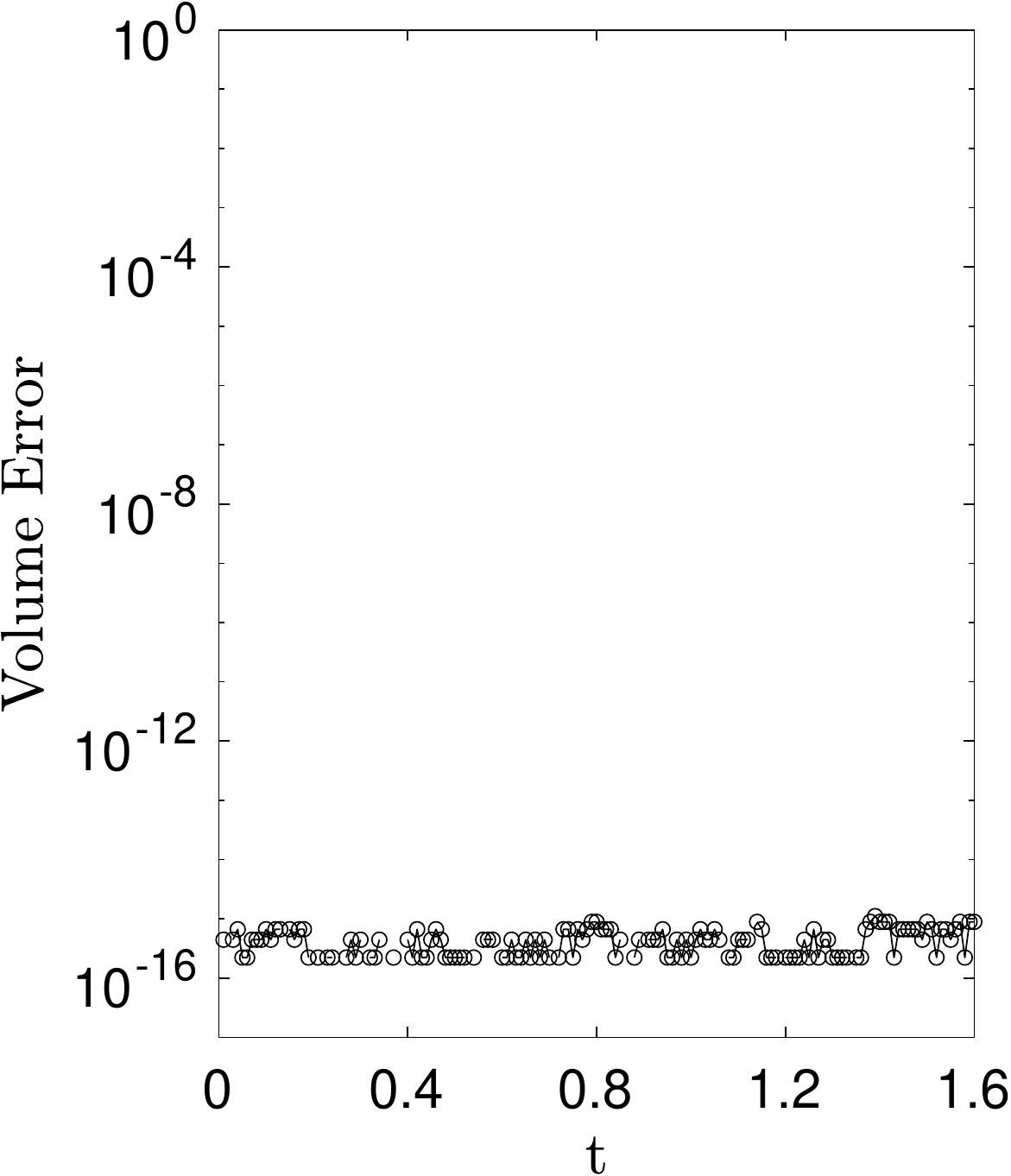}}
\end{center}
\end{minipage}
\end{minipage}
\caption{The evolution of the $\tanh$ from $0 \leq t \leq 1.6$ with the flow field containing cosine multiplier. The dash line is the function at $t=0$. The left most figure is the profile at time $t=0.8$ ($t =T_\textnormal{fin}/2$) and the middle figure is at time $t=1.6$ ($t =T_\textnormal{fin}$). The right most figure shows the absolute error between the volume at time $t$ with the one computed at initial time. The 1D blob does seem to coincide back with its original shape at the final time. Although unnecessary for this purpose, the simulation trajectory is precomputed with \texttt{solvebvp} in Chebfun.}
\label{fig:numex2}
\end{figure}

This test is similar to the Experiment 1 with the flow field given by $u(x,t) = c (2 \pi / T_\textnormal{fin}) x  \cos(2 \pi t / T_\textnormal{fin})$. The cosine multiplier function has the effect of decreasing the magnitude of the vector field until time $t=T_\textnormal{fin}/2$ and then powering back up in a reversible way such that the evolving function $f$ should coincide with the original shape at $t=0$ at time $T$. This test can give us insights into whether the method can maintain the shape and volume or whether some smearing effect is happening. \cref{fig:numex2} shows that the volume error is near machine zero, and the original shape is maintained at the final time with no smearing. Although we can compute $x(t)$ analytically, we can instead use an ODE solver to solve it. For every endpoint/vertex of the patches, we precompute their trajectory offline. For example, using the Chebfun ODE solver for the left endpoint of the patch with index 1 ($\Omega_1$), we can use the commands

\begin{verbatim}

N = chebop(Tinit, Tfin);
N.op = @(t,x) diff(x)-x*c*(2*pi/Tfin)*cos(2*pi*t/Tfin);
N.lbc = f.leafArray{1}.domain(1);
odesol = solvebvp(N,0) or odesol = N\0.

\end{verbatim}
The \texttt{odesol} is a function handle that defines the path line or trajectory of the left endpoint of the patch $\Omega_1$, e.g., \texttt{odesol(0.5)} will compute its position at time $t=0.5$. The \texttt{solvebvp} is essentially solving the initial value problem in ``space-time" where $t$ is treated the same way as the spatial variable. As the results demonstrate, the method maintains accuracy down to machine precision in both space and time. 

\subsection{Experiment 3 (pure strain flow in 2D)}
\begin{figure}[htbp]
\begin{minipage}{\linewidth}
\begin{minipage}{0.475\linewidth}
\scalebox{0.375}{\includegraphics{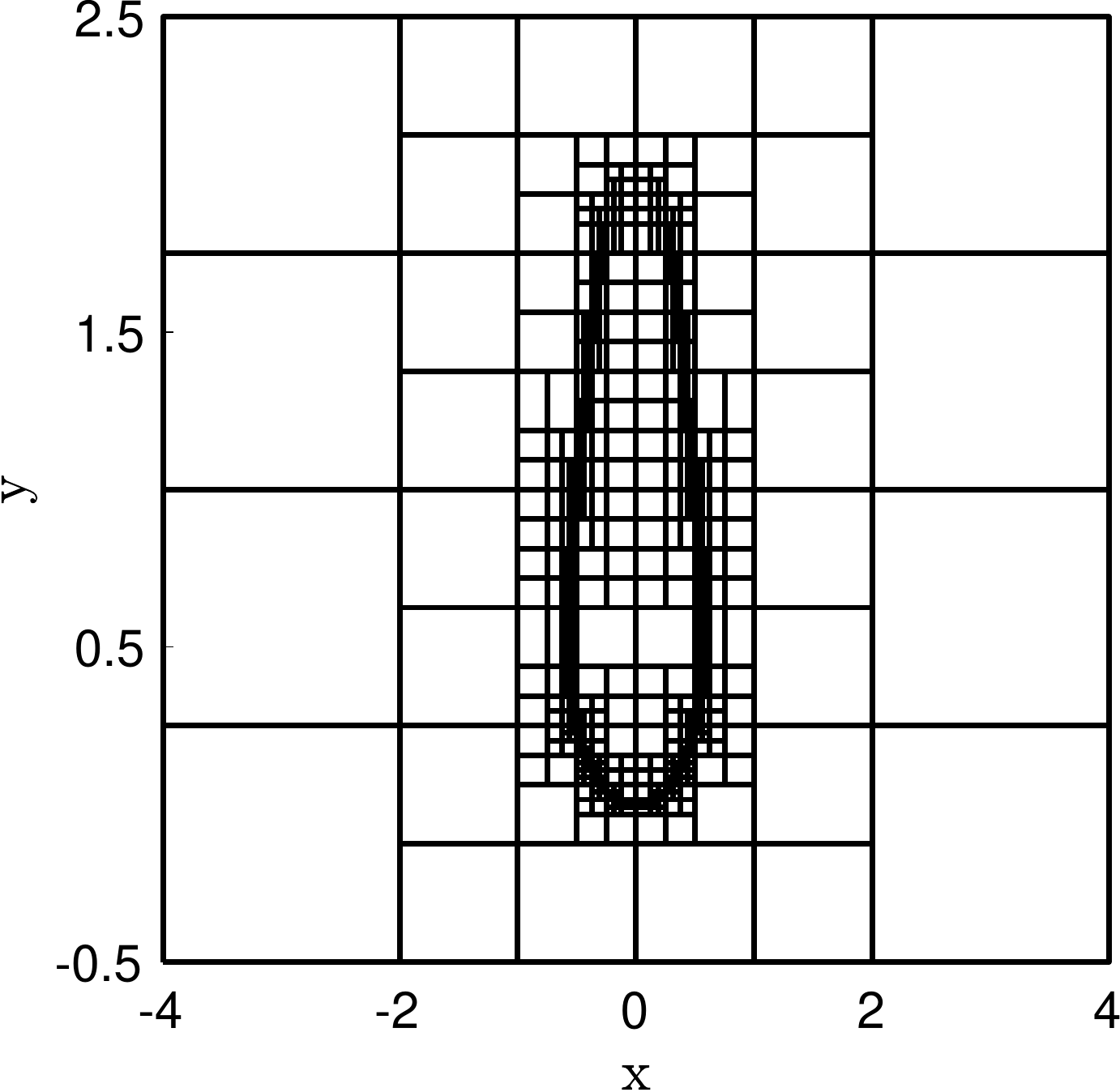}}
\end{minipage}
\hspace{10pt}
\begin{minipage}{0.475\linewidth}
\scalebox{0.375}{\includegraphics{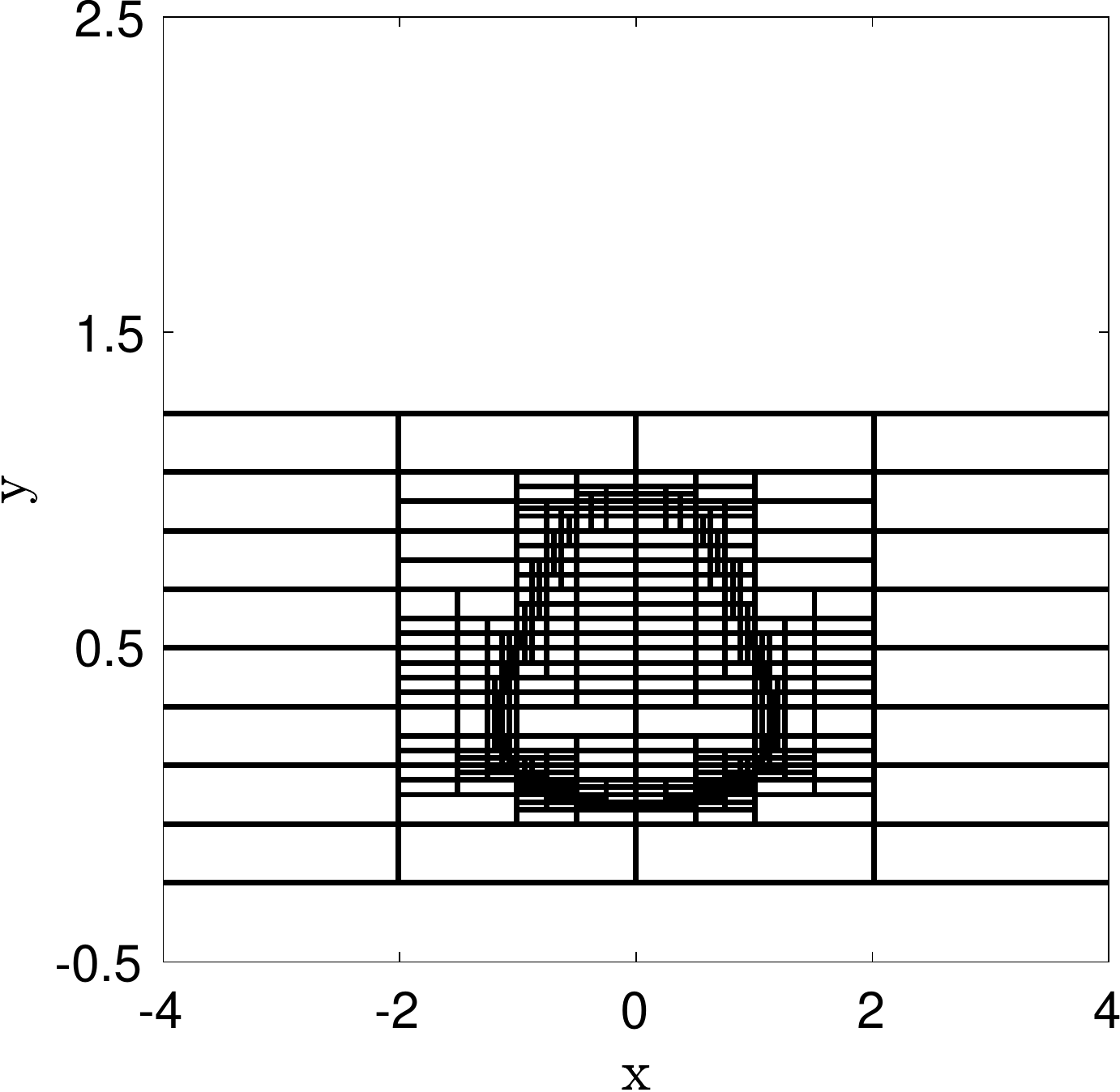}}
\end{minipage}
\end{minipage}
\caption{Left: the pear-shaped form dynamics under the pure strain flow. The partition of unity patches distribution follows the profile of the pear. Patches with bigger sizes are away outward from the interface. Right: since we let the corners of the domain to freely move, the patches at $t=0.8$ are squeezed towards $y=0$.}
\label{fig:numex3}
\end{figure}

\begin{figure}[htbp]
\begin{minipage}{\linewidth}
\begin{minipage}{0.49\linewidth}
\centering $t=0$
\scalebox{0.32}{\includegraphics{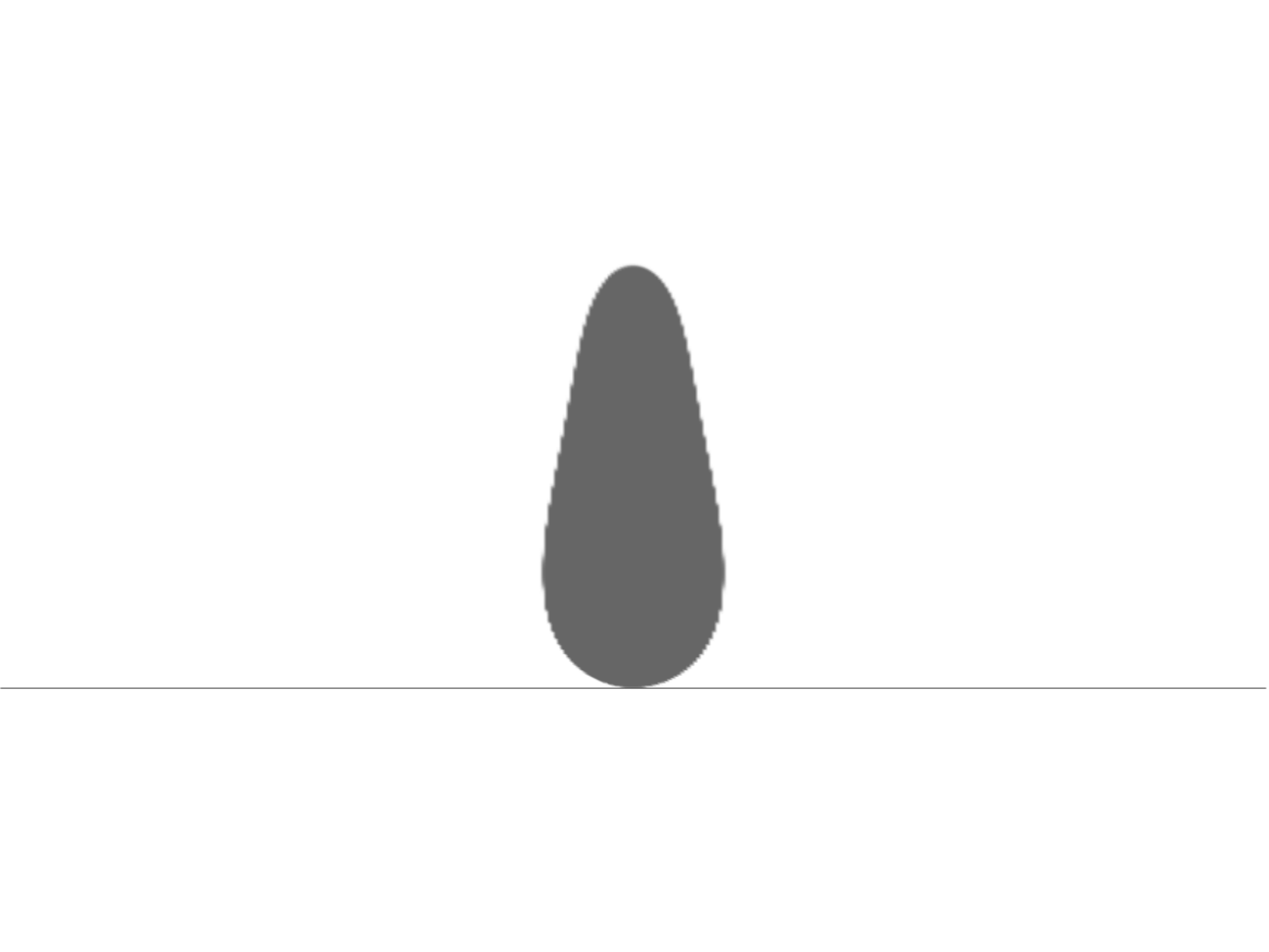}}
\end{minipage}
\hspace{3pt}
\begin{minipage}{0.49\linewidth}
\centering $t=0.8$
\scalebox{0.32}{\includegraphics{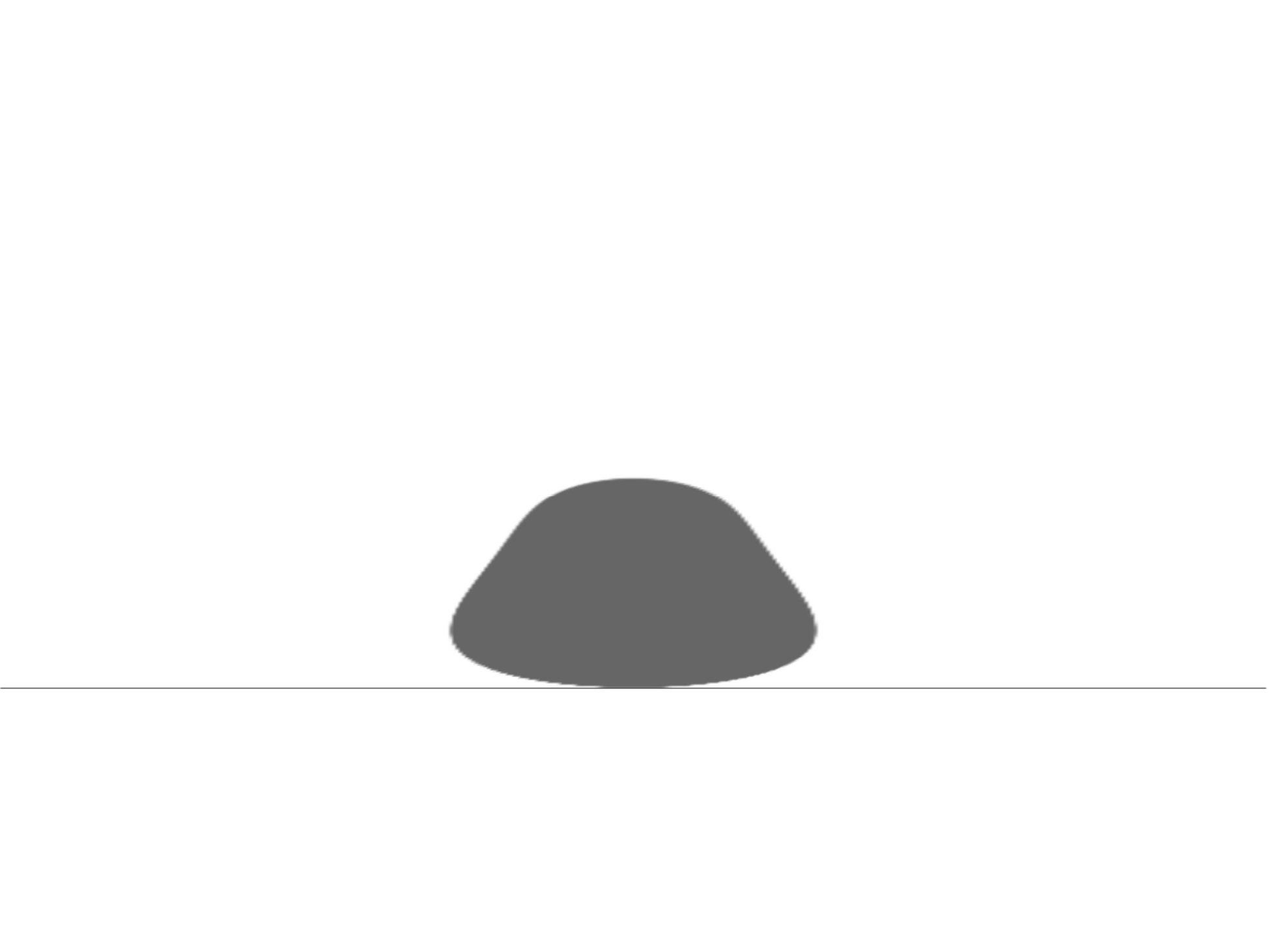}}
\end{minipage}
\end{minipage}
\begin{minipage}{\linewidth}
\begin{minipage}{0.49\linewidth}
\centering $t=1.6$
\scalebox{0.32}{\includegraphics{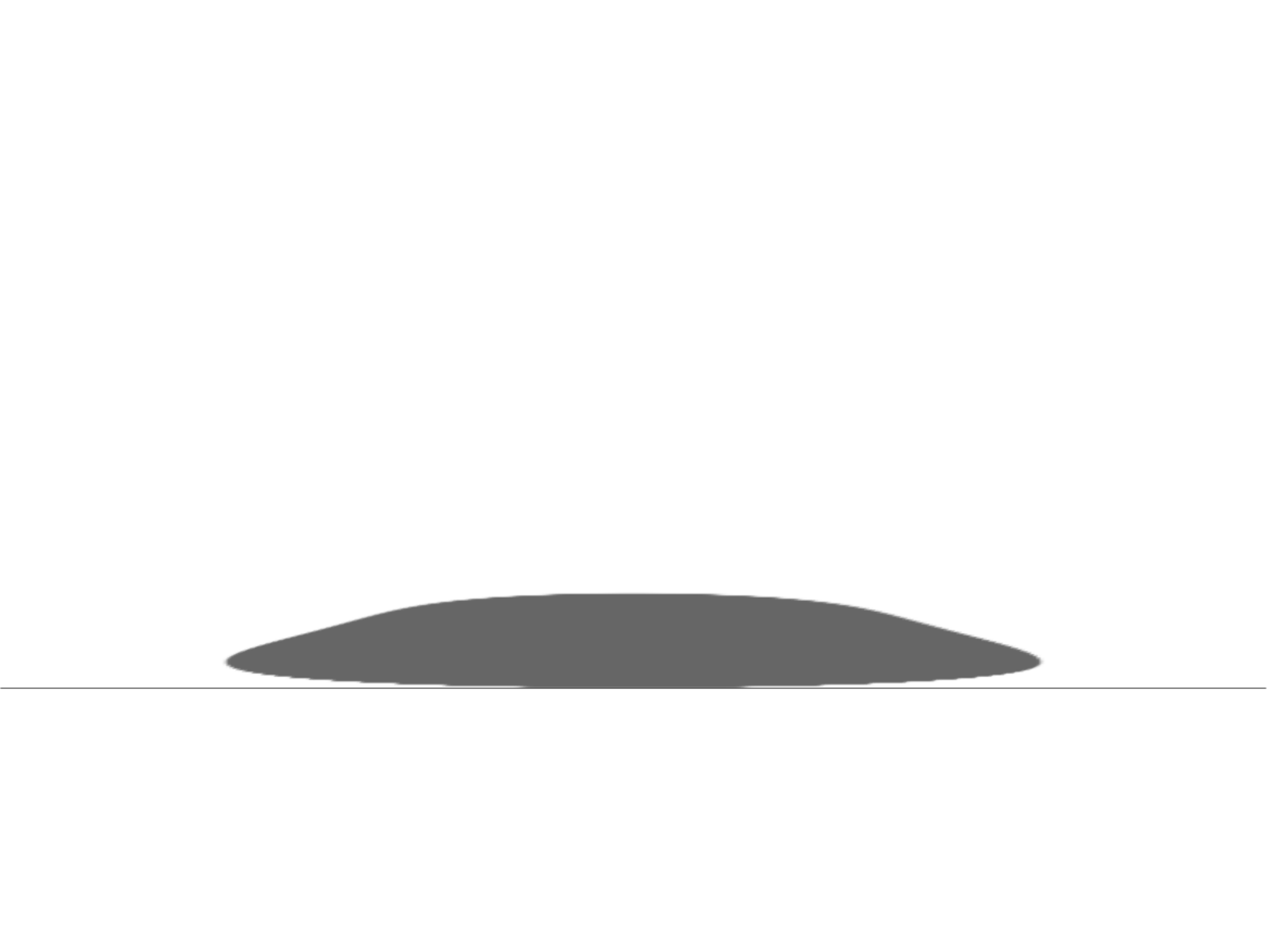}}
\end{minipage}
\hspace{25pt}
\begin{minipage}{0.4\linewidth}
\scalebox{0.375}{\includegraphics{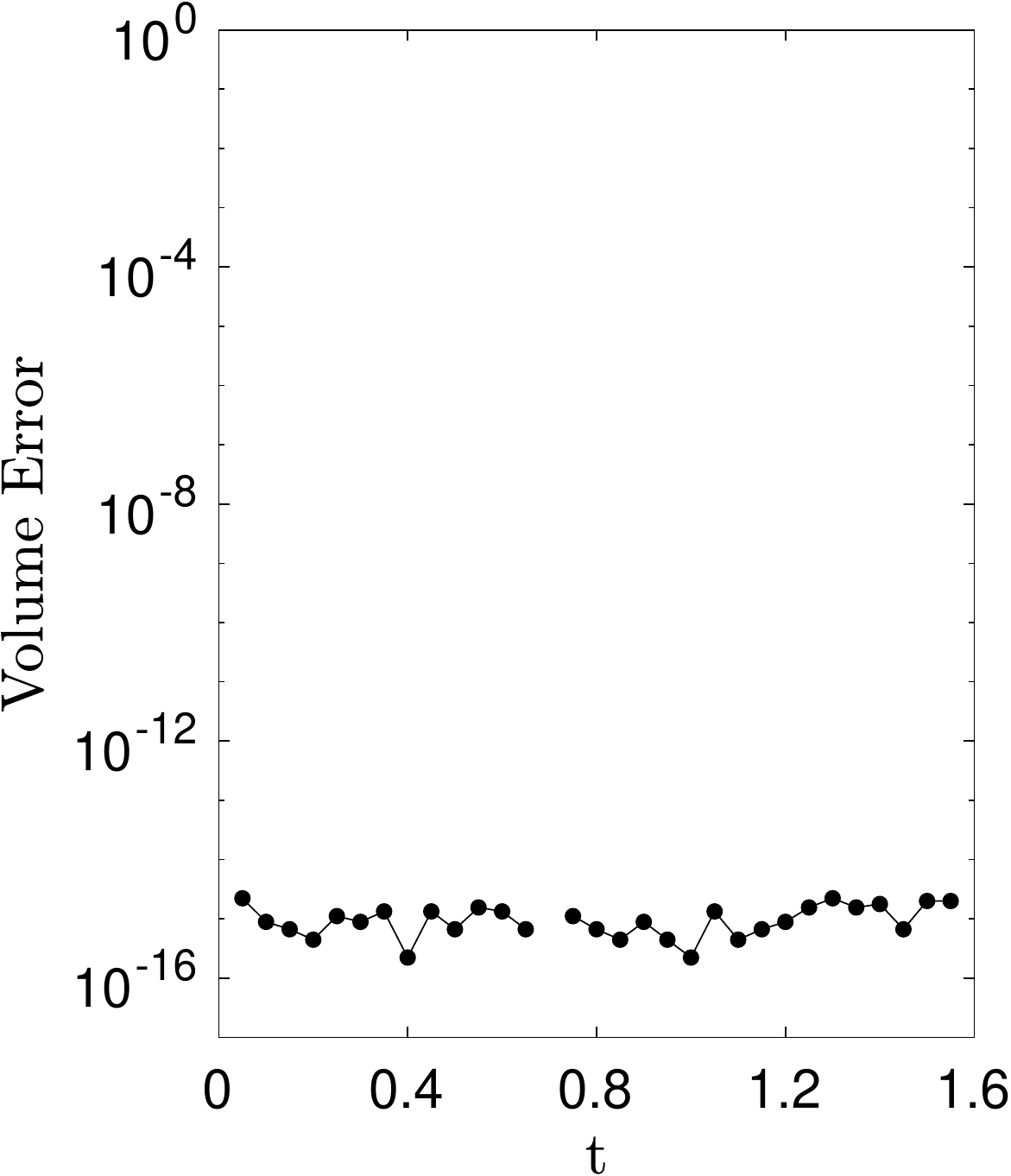}}
\end{minipage}
\end{minipage}
\caption{Top left: initial state of the pear shape at $t=0$. The pear condition at $t=0.8$ ($t =T_\textnormal{fin}/2$) (top right) and at the final time $t=1.6$ ($t =T_\textnormal{fin}$) (bottom left). }
\label{fig:numex3b}
\end{figure}

For this experiment, we use a 2D function with a pear-shaped form given by
\begin{displaymath}
f(x,y) = \tfrac{1}{2}\tanh\left(100\left(\left(1+\left(y-1\right)^3\right)\left(1-\left(y-1\right)\right)-4x^2\right)\right)+\tfrac{1}{2}
\end{displaymath}
on the domain $[-4,4] \times [-0.5,2.5]$. This pear shape, shown in \cref{fig:numex3b}, is a class of piriform curves commonly used to approximate droplet shape. The flow field is of type pure strain $\underline{u} = (x,-y)$. \cref{fig:numex3} shows the initial shape and its initial distribution of patches. For this simulation, we let the vertices of the domain move along with the flow. Under this flow, we expect the pear to be squeezed/deformed towards the line $y=0$, and the patches become elongated in the horizontal direction at the final time. \cref{fig:numex3b} shows the volume accuracy is still down to machine precision throughout the simulation.

\subsection{Experiment 4 (circular motion without rotation)}
In this experiment, we perform a circular motion of the medium \emph{without} rotation. The object (material/medium 1)  is a unit disk initially centered at $(0,2)$. The function to represent it is 
\begin{displaymath}
f(x,y) = \tfrac{1}{2} \tanh(50(1-(x^2+(y-2)^2))) + \tfrac{1}{2}.
\end{displaymath}
With the steep gradient slope of $50$, the area of the disk computed with \texttt{sum(f)} is already close to $\pi$ up to $14$ digits. We then let the disk move in circular motion such that the distance of its center of mass to the origin is $2$ throughout the simulation. Additionally, the vertices of the domain are kept fixed. Hence, the patches are constrained to move only inside the bounding box. 

\begin{figure}[htbp]
\begin{minipage}{\linewidth}
\scalebox{0.35}{\includegraphics{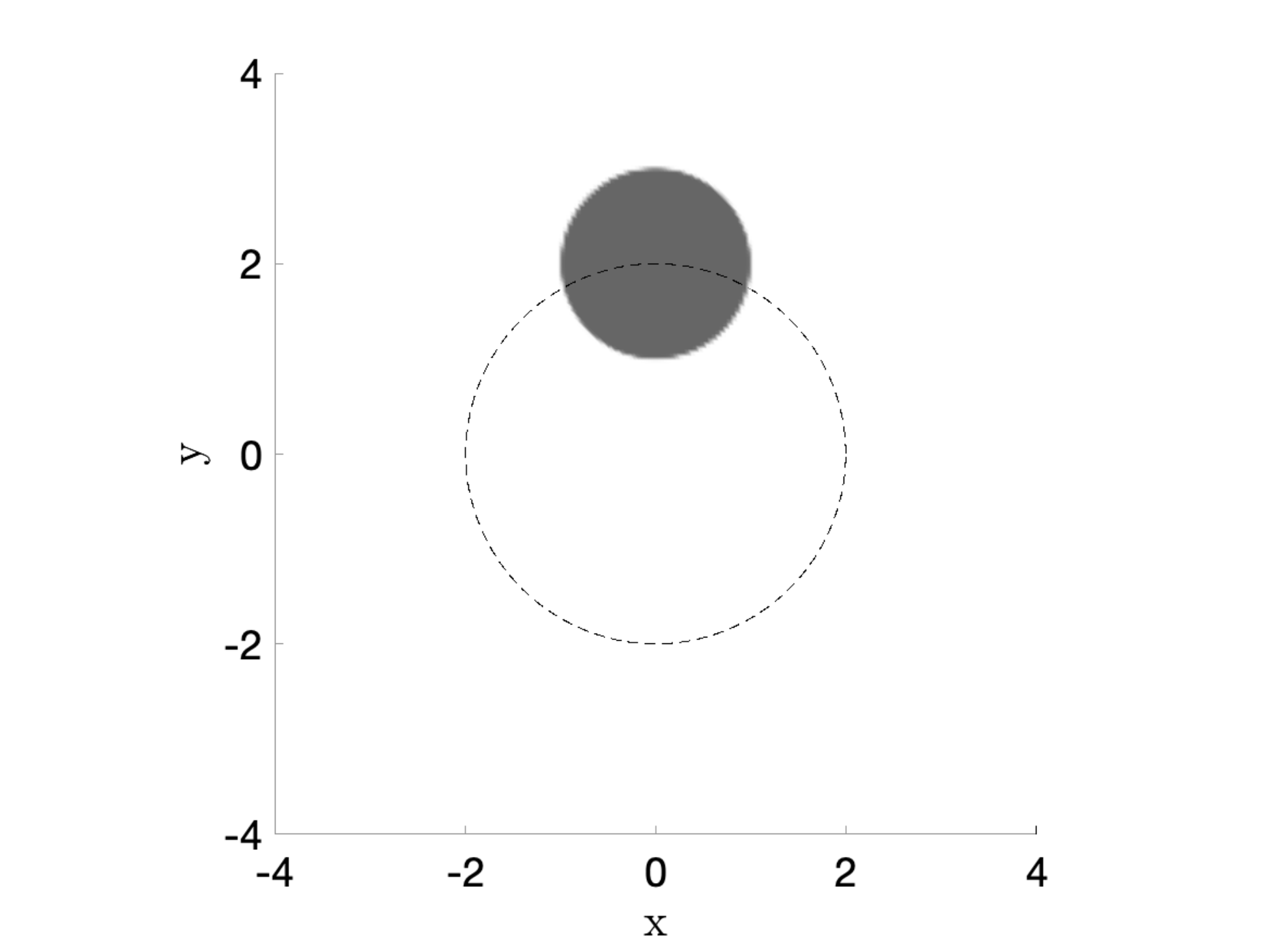}}
\hspace{-30pt}
\scalebox{0.35}{\includegraphics{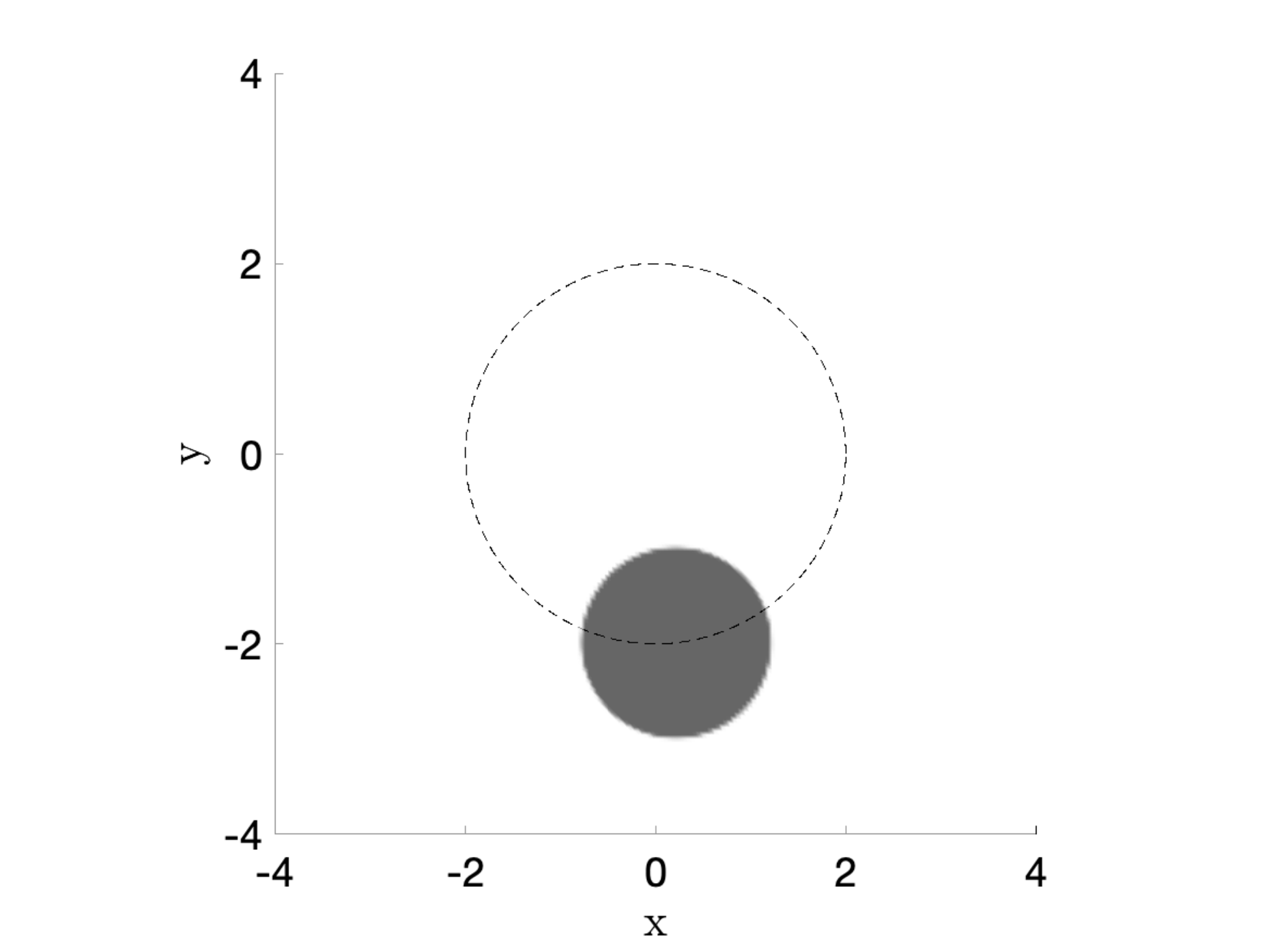}}
\end{minipage}
\begin{minipage}{\linewidth}
\vspace{10pt}
\begin{center}
\scalebox{0.375}{\includegraphics{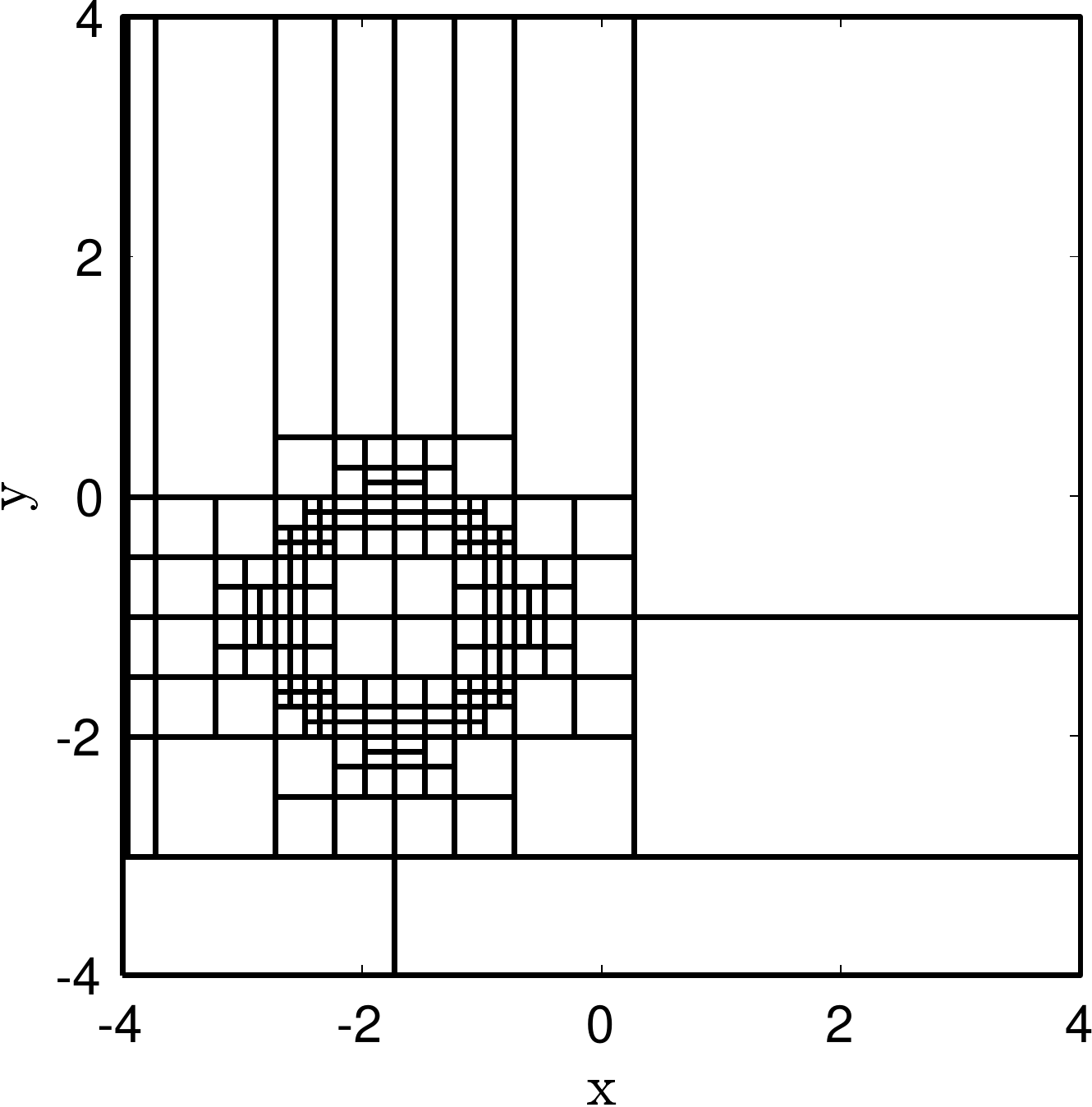}}
\hspace{25pt}
\scalebox{0.375}{\includegraphics{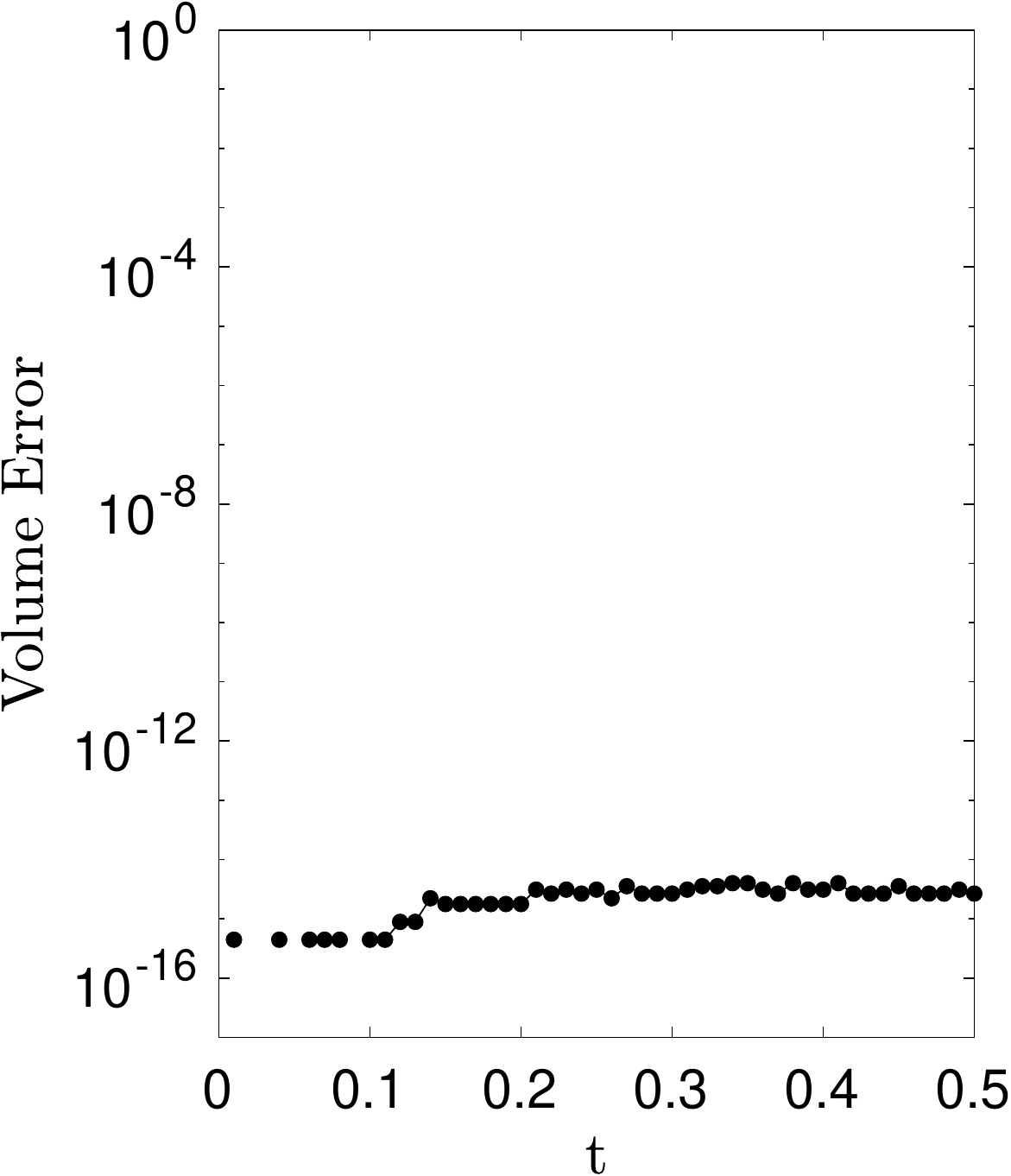}}
\end{center}
\end{minipage}
\caption{The unit disk in circular motion. The distance of its center of mass to the origin is $2$ throughout the simulation. The center of mass stays on the orbit (i.e., a circle of radius 2) and the volume is conserved.}
\label{fig:numex4}
\end{figure}

The top two figures of \cref{{fig:numex4}} show the position of the unit disk at the initial time and the half-circular way of the motion. The patches are shown when the position of the disk is somewhere in the third quadrant. Since the vertices of the bounding box are kept fixed, some patches are stretched to maintain good domain coverage. The volume accuracy can be seen in the bottom right plot of the \cref{fig:numex4}. Although the location of the center of mass is available analytically during the motion, one can also check it numerically. Note that the data structure contains all information about Chebsyhev grids and all necessary weights and coefficients in each patch. Hence, one can use them to compute the center of mass or centroids in each patch or other type of ``moments". 

Even though the method is not yet intended for interface reconstruction, we can utilize the interpolant, its derivatives, and integrals to compute the volume fractions, centroids, center of mass, and gradients needed for a very rudimentary reconstruction. The formula provided in \cref{tab:moments} can be utilized.
\begin{table}[htbp]
\footnotesize
\caption{Volume fraction, centroid, center of mass, unit gradient vector}\label{tab:moments}
\begin{center}
 \begin{tabular}{|c|c|} \hline
  Volume fraction & $v_f(\Omega_i) = \frac{\int f d\Omega_i}{\textnormal{max}(f)_{\Omega_i} vol(\Omega_i)}$ \\
  &\\
  Centroid & $\underline{x}_c = \left( \frac{\int x f d\Omega_i}{vol(\Omega_i)},\frac{\int y f d\Omega_i}{vol(\Omega_i)} \right)$ \\ 
  &\\
  Center of mass &  $\underline{x}_m = \left( \frac{\int x f d\Omega}{vol(\Omega)},\frac{\int y f d\Omega}{vol(\Omega)} \right)$ \\ 
  &\\
  Unit gradient &  $\underline{\widehat{g}} = \frac{\nabla f}{|\nabla f|}$ \\ 
  \hline
 \end{tabular}
\end{center}
\end{table}
The use of volume fraction for basic interface reconstructions is in the spirit of the volume-of-fluid (VOF) method \cite{hirt_volume_1981} and its variants. Moreover, the use of a root-finding strategy in the direction of unit gradients to find a point in the neighborhood of the interface is common in the Moment-of-Fluid methods \cite{jemison_coupled_2013, chen_predicted-newtons_2019,maric_iterative_2021} and its vast literature.

\begin{figure}[htbp]
\begin{minipage}{\linewidth}
\scalebox{0.39}{\includegraphics{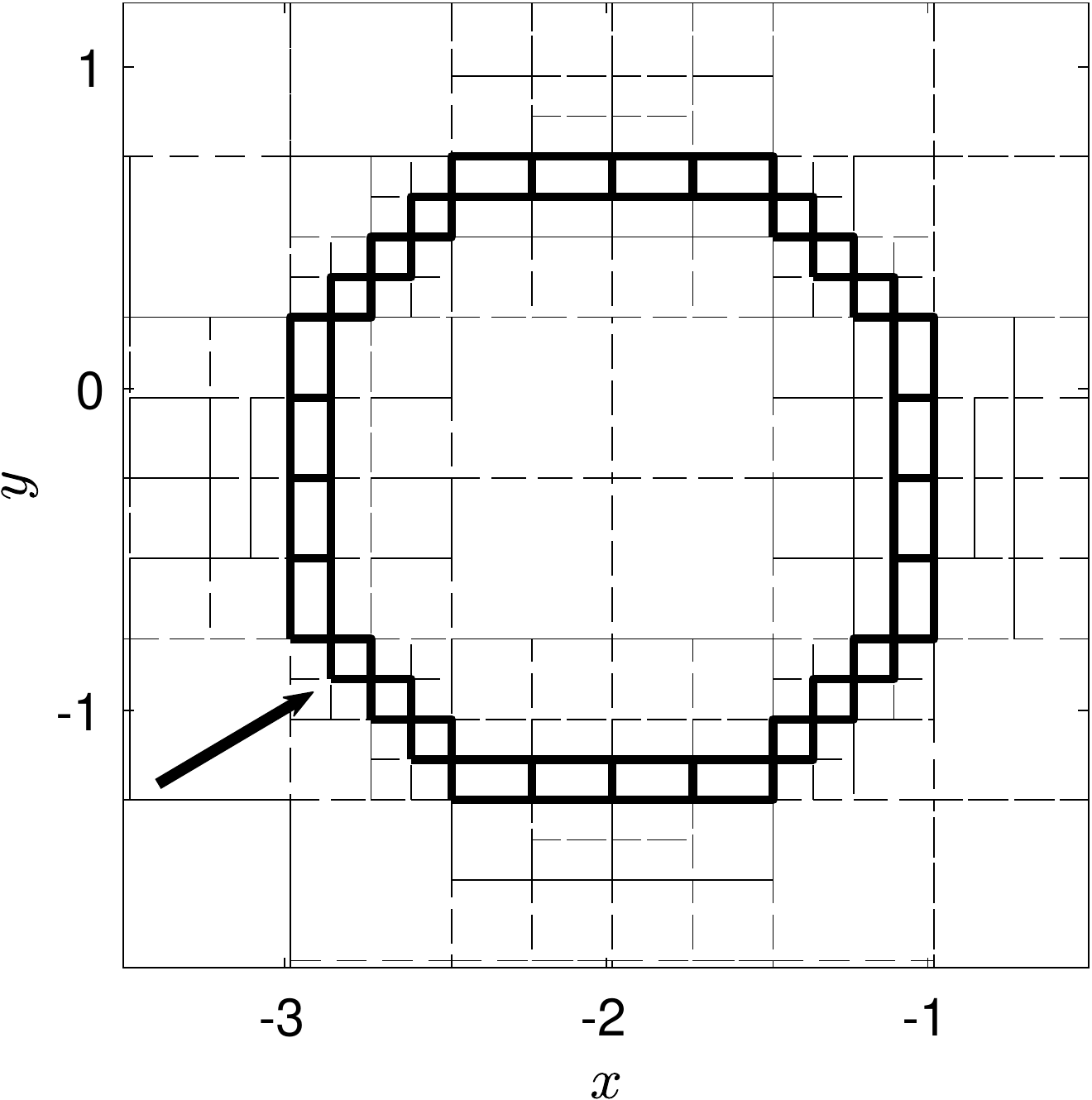}}
\scalebox{0.39}{\includegraphics{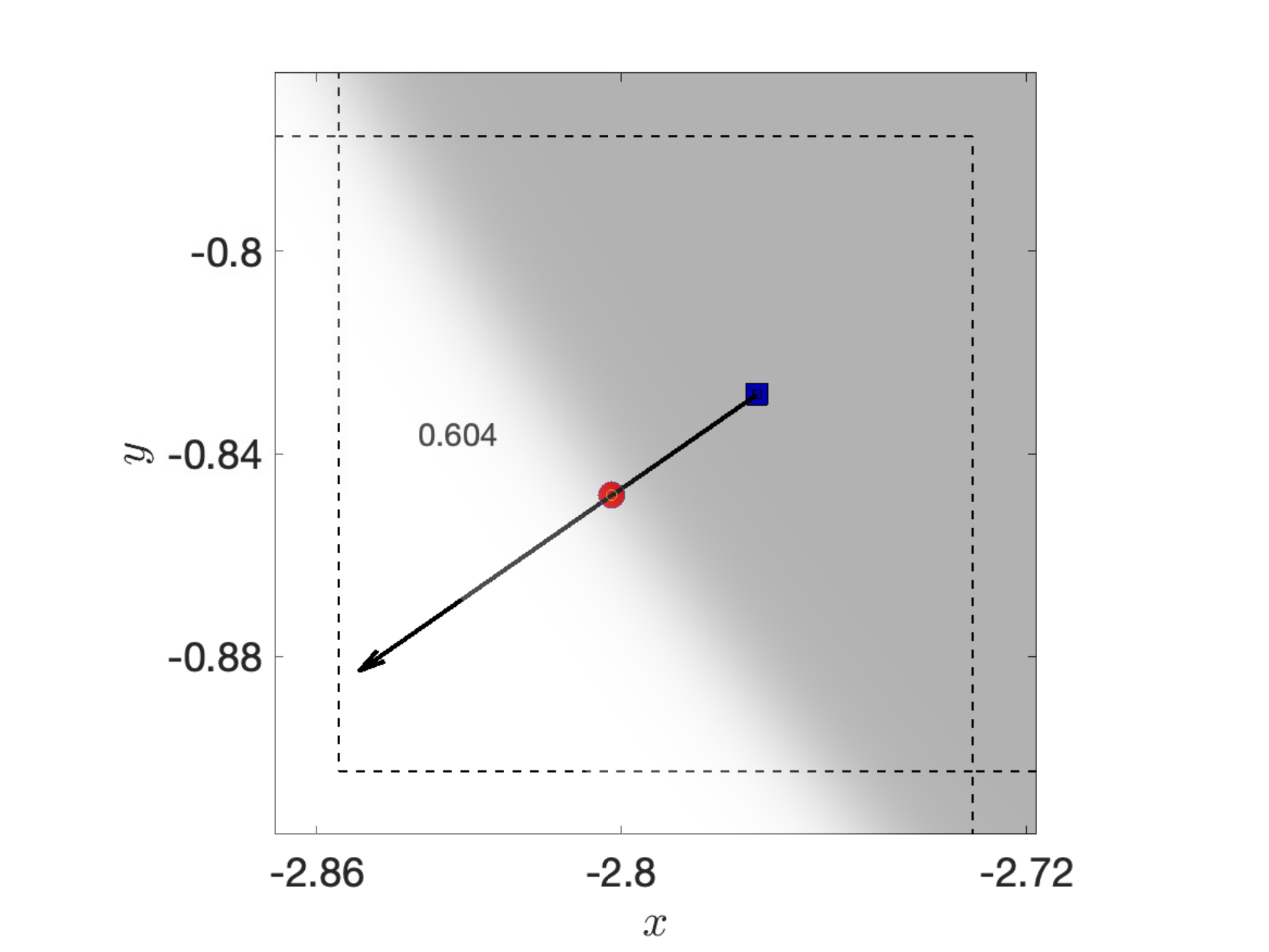}}
\end{minipage}
\begin{minipage}{\linewidth}
\centering
\scalebox{0.45}{\includegraphics{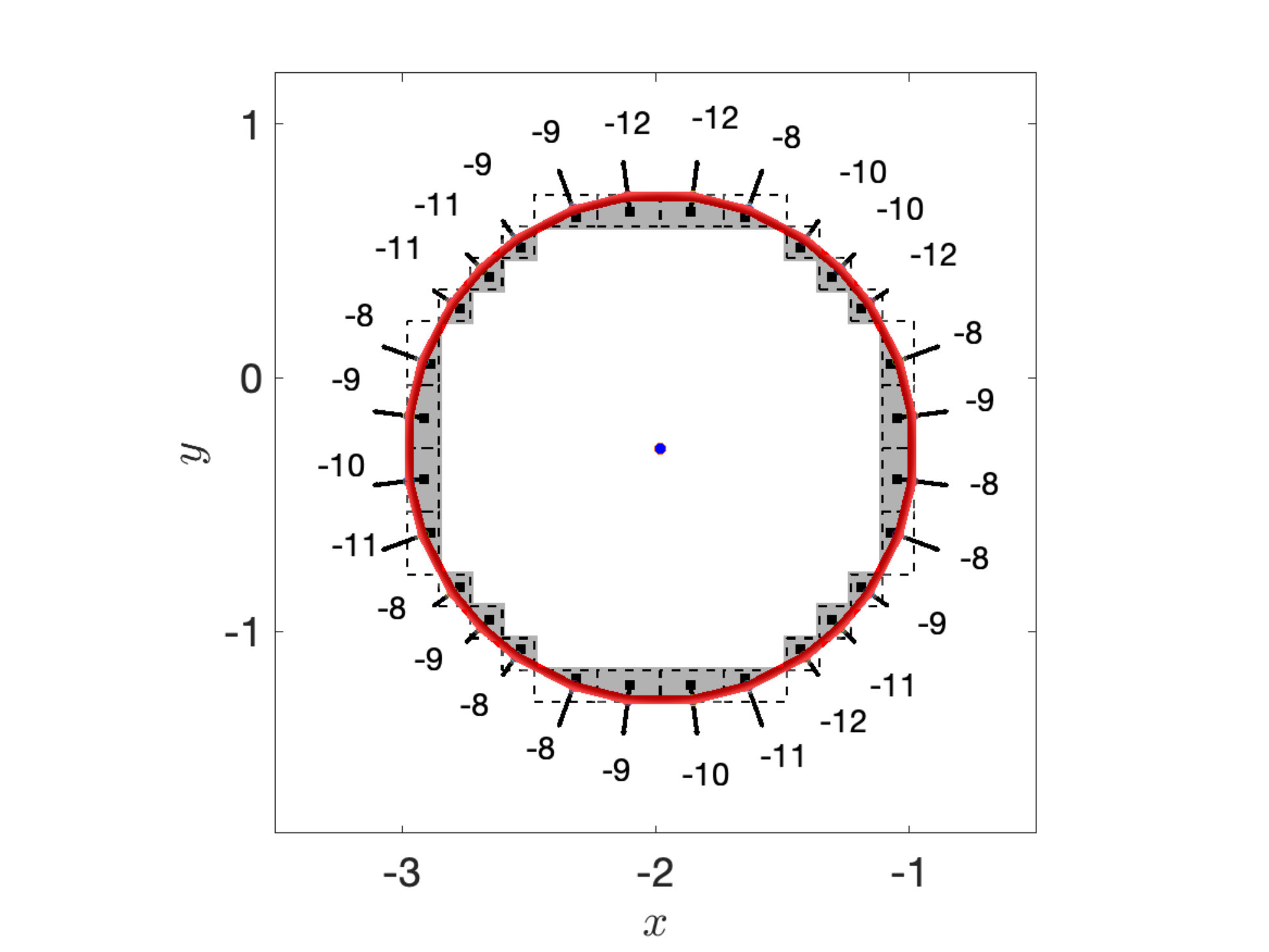}}
\end{minipage}
\caption{Top left: a zoom-in view of the the disk somewhere around the third quadrant from \cref{fig:numex4}. Boxes/zones with volume fractions range from 5\% to 95\% are in bold. Top right: a zoom-in view of one of the patches pointed by the arrow on the top left figure. The patch has a volume fraction of about 60.4\%. The point on the interface is obtained using root-finding method applied to the parametric equation of a line starting from the centroid in the direction of the unit gradient outward. Bottom: points on the interface are then connected to recover an approximate curve (in red) of the interface. The error of the curvature, i.e. the $\log_{10}(|\kappa - 1|)$, with respect to the exact curvature at the boundary of a unit disk is also shown. The center of mass location is accurate to 14 digits.}
\label{fig:numex4curvekappa}
\end{figure}

We first compute volume fractions of all patches and keep the ones between, say, 5\% and 95\%. The patches containing volume fractions between $0$ and $1$ are shared between two mediums. The top left figure of \cref{fig:numex4curvekappa} shows patches/zones in bold where the volume fractions range between 5\% and 95\%. We can compute the centroids and the unit gradient vectors (at the centroids) for those patches. A zoom-in view of one of the patches, pointed by the arrow, is shown in the top right figure of \cref{fig:numex4curvekappa}. In that particular patch, the volume fraction is about 60.4\%.

To find an approximate interface point in each patch, we can search for a point away from the centroid in the direction of the gradient vector using the following parametric equation.
\begin{displaymath}
\underline{x}(\tau) = \underline{x}_c + \tau \textnormal{diam}(\Omega_c)\underline{\widehat{g}}_c,
\end{displaymath}
where $\textnormal{diam}(\Omega_c)$ is the diameter of the patch/zone where the centroid $x_c$ is located. We can then use the APU interpolant $f$ and root-finding method for finding parameter $\tau$ such that $f(\underline{x}_\textnormal{ifc}) = 1/2$. Essentially, we find points where the values of the function $f$ are $1/2$ (a halfway value of the two mediums). For VOF-based flow solvers like ours, it is common to plot an iso-contour of $0.5$ of volume fractions. On a particular patch, as in the top right figure of \cref{fig:numex4curvekappa}, the point on the interface is colored in red. After those points are found, we can translate and scale their distances relative to the center of mass, transform them into polar coordinates, and sort their angles. One can then plot the points in a counter-clockwise way to get the rough approximate curve of the interface, as shown at the bottom of \cref{fig:numex4curvekappa}.

An important geometrical quantity in simulations of multiphase flows is the interface curvature, which directly determines the surface tension force along the interface at any given point. Next, the values of interface curvature at those points can be computed with (see \cite{goldman_curvature_2005})
\begin{displaymath}
\kappa = \frac{|f_y^2f_{xx} - 2f_xf_yf_{xy}+f_x^2f_{yy} |}{(f_x^2 + f_y^2)^{3/2}}
\end{displaymath}
In this case, the exact curvature for all points at the boundary of the unit disk is $\kappa = 1$. In our computations, the curvature at those points is between $8$ to $12$ digits accuracy. \cref{fig:numex4curvekappa} shows the $\log$ of absolute error of the curvature at each computed interface point. 

\subsection{Experiment 5 (rotation of multibodies)}
\begin{figure}[htbp]
\begin{minipage}{\linewidth}
\hspace{-10pt}
\scalebox{0.38}{\includegraphics{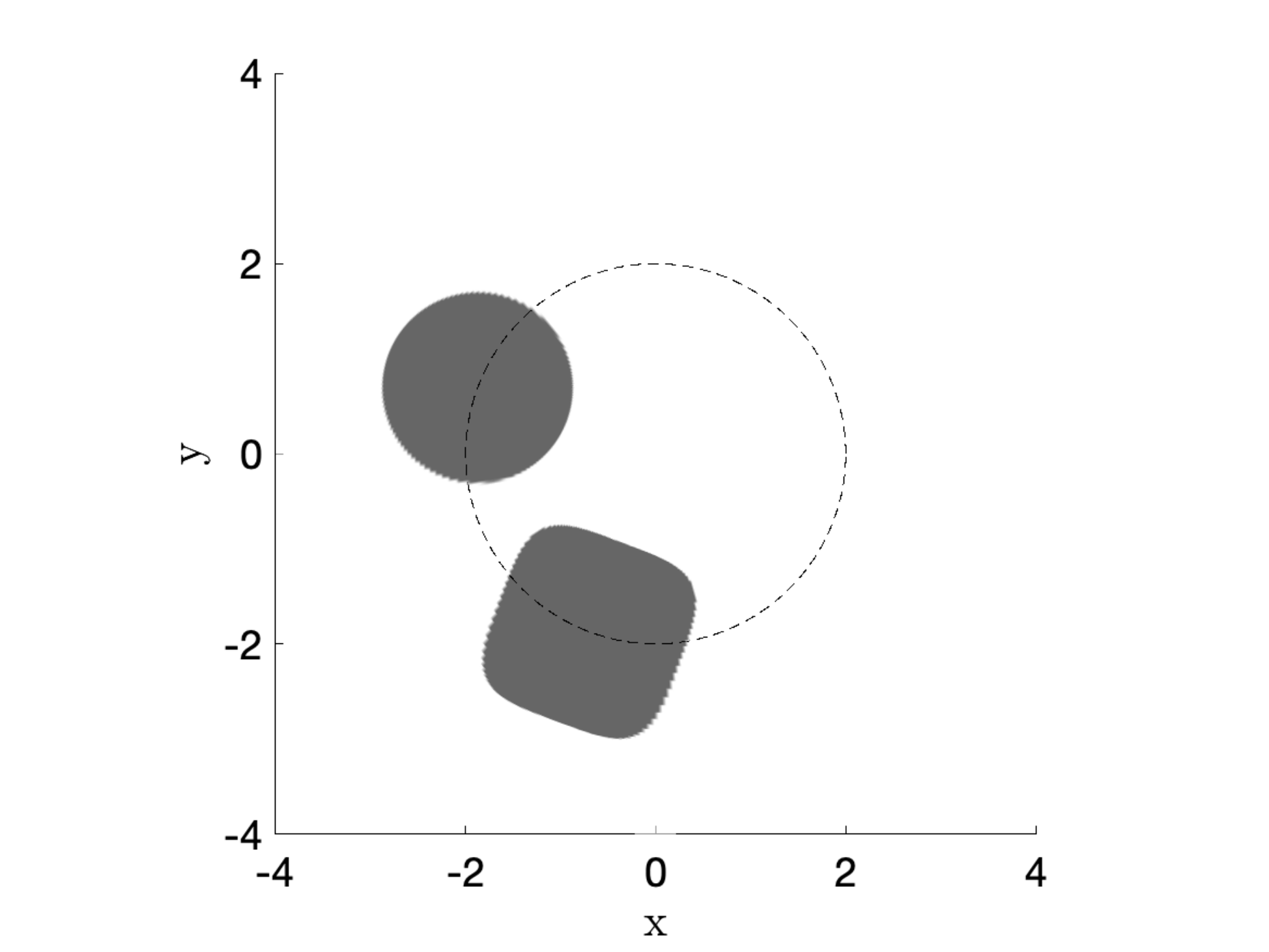}}
\hspace{-55pt}
\scalebox{0.38}{\includegraphics{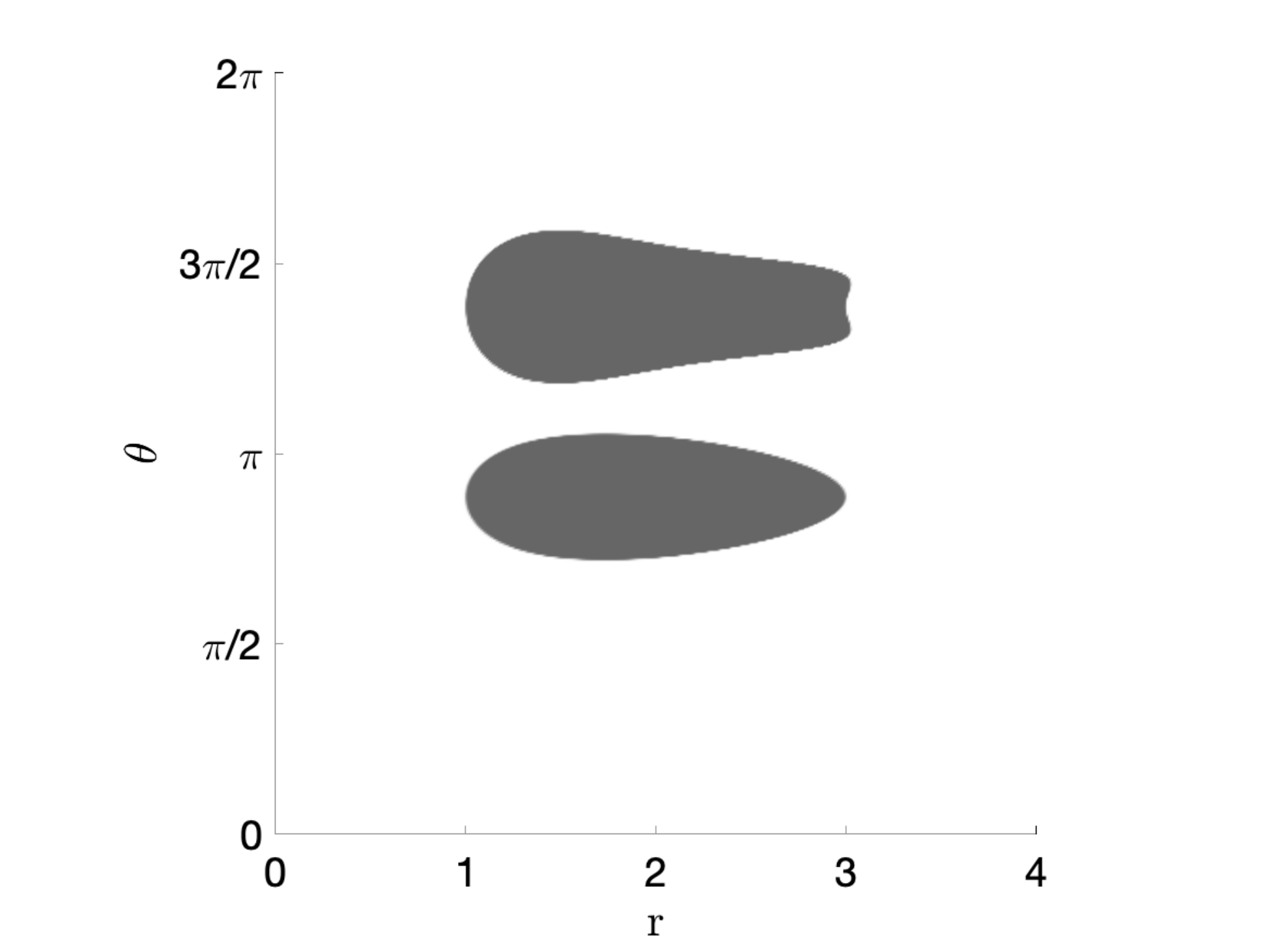}}
\end{minipage}
\begin{minipage}{\linewidth}
\vspace{10pt}
\begin{minipage}{0.48\linewidth}
\centering
\scalebox{0.4}{\includegraphics{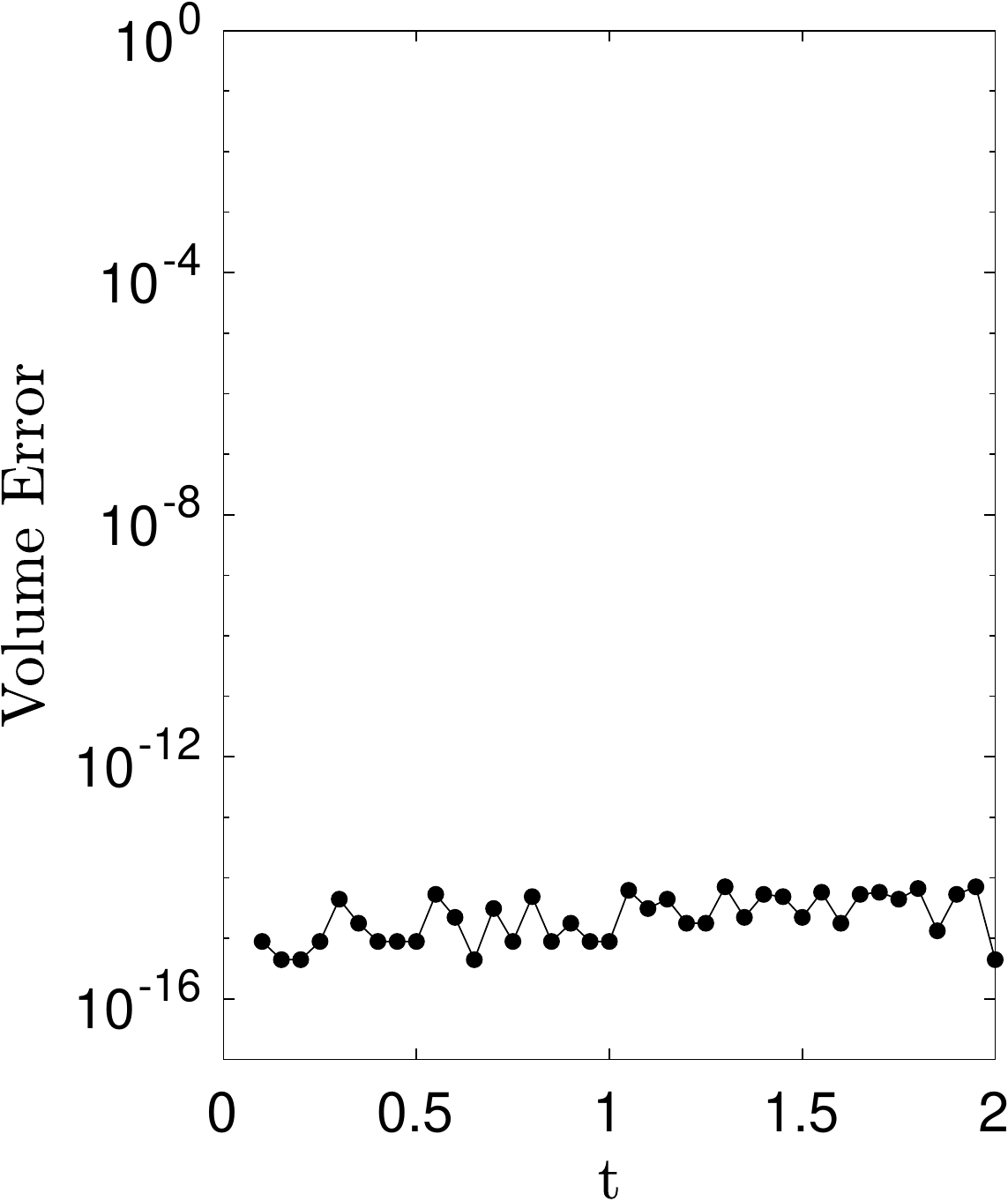}}
\end{minipage}
\begin{minipage}{0.48\linewidth}
\centering
\scalebox{0.4}{\includegraphics{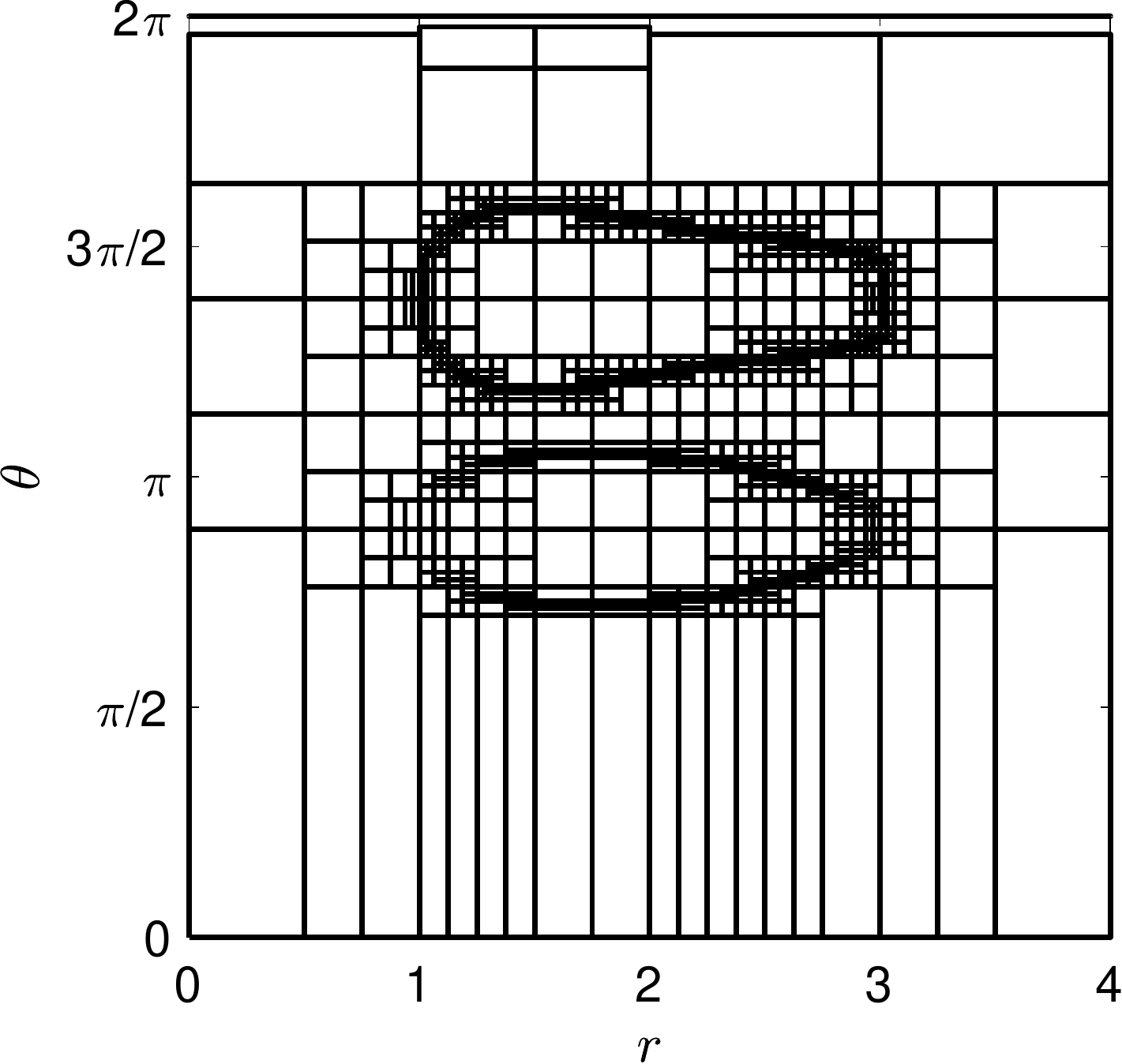}}
\end{minipage}
\end{minipage}
\caption{Rotations of a circle and a squircle under unit pure rotational vector field. In polar coordinates, the shapes are an egg-like and a pear-like geometry. The top left figure shows the position of both mediums at $t=2$. The top right figure shows their positions at $t=2$ in polar coordinates, where the motion becomes a simple translation. The patches distributions at $t=2$ shows level of details around the interfaces.}
\label{fig:numex5}
\end{figure}

For Experiment 5, we are rotating a unit disk and a unit squircle together (both centered along the circle of radius $2$) under a pure rotation unit vector field $\underline{u}=(-y/r,x/r)$, with $r=\sqrt{x^2 + y^2}$ from time $t \in [0,2]$. The two-medium function is represented in $r$ and $\theta$ at time $t=0$ with a slope of $100$ as
\begin{align*}
f(r,\theta) = & 1 + \tfrac{1}{2}\tanh(100(1-((r\cos(\theta+\pi/4))^2+(r\sin(\theta+\pi/4)-2)^2))) + \\
&\tfrac{1}{2}\tanh(100(1-((r\cos(\theta-\pi/4))^4+(r\sin(\theta-\pi/4)-2)^4)))+1
\end{align*}
If the computation is done in the cartesian coordinates, then rectangular patches are going to break up due to the rotations. However, in the polar domain, the rotation becomes a translation in the angular direction as $\underline{u}(r,\theta)=(0,1)$. \cref{fig:numex5} also shows that the circle and squircle, once transformed into the polar domain, become an egg-like and pear-like geometry. 

The circle and the squircle are positioned a bit closer together on purpose. This tests whether the APU constructors can correctly handle and adaptively refine regions between the two interfaces. The simulation shows no degradations of the total volume and no changes in shapes. When the pear and the egg are close to $\theta = 2\pi$, one can shift the domain upward to keep the rotation since the interpolant is initially constructed for non-periodic domains. 

\subsection{Experiment 6 (angular deformation with cosine multiplier)}

\begin{figure}[htbp]
\begin{minipage}{\linewidth}
\begin{minipage}{0.49\linewidth}
\begin{center}
$t=0.25$ \\
\scalebox{0.32}{\includegraphics{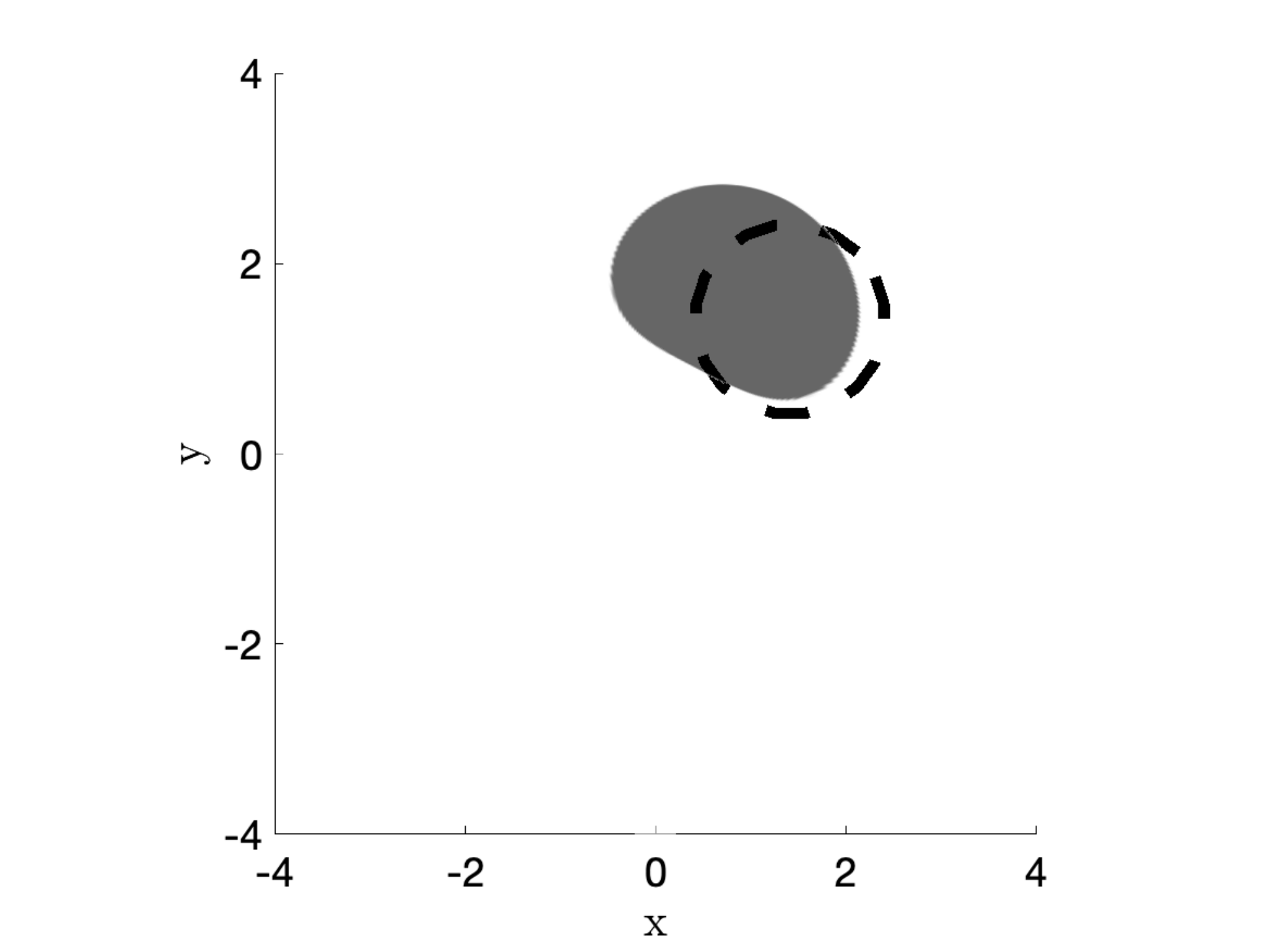}}
\end{center}
\end{minipage}
\begin{minipage}{0.49\linewidth}
\begin{center}
$t=0.5$ \\
\scalebox{0.32}{\includegraphics{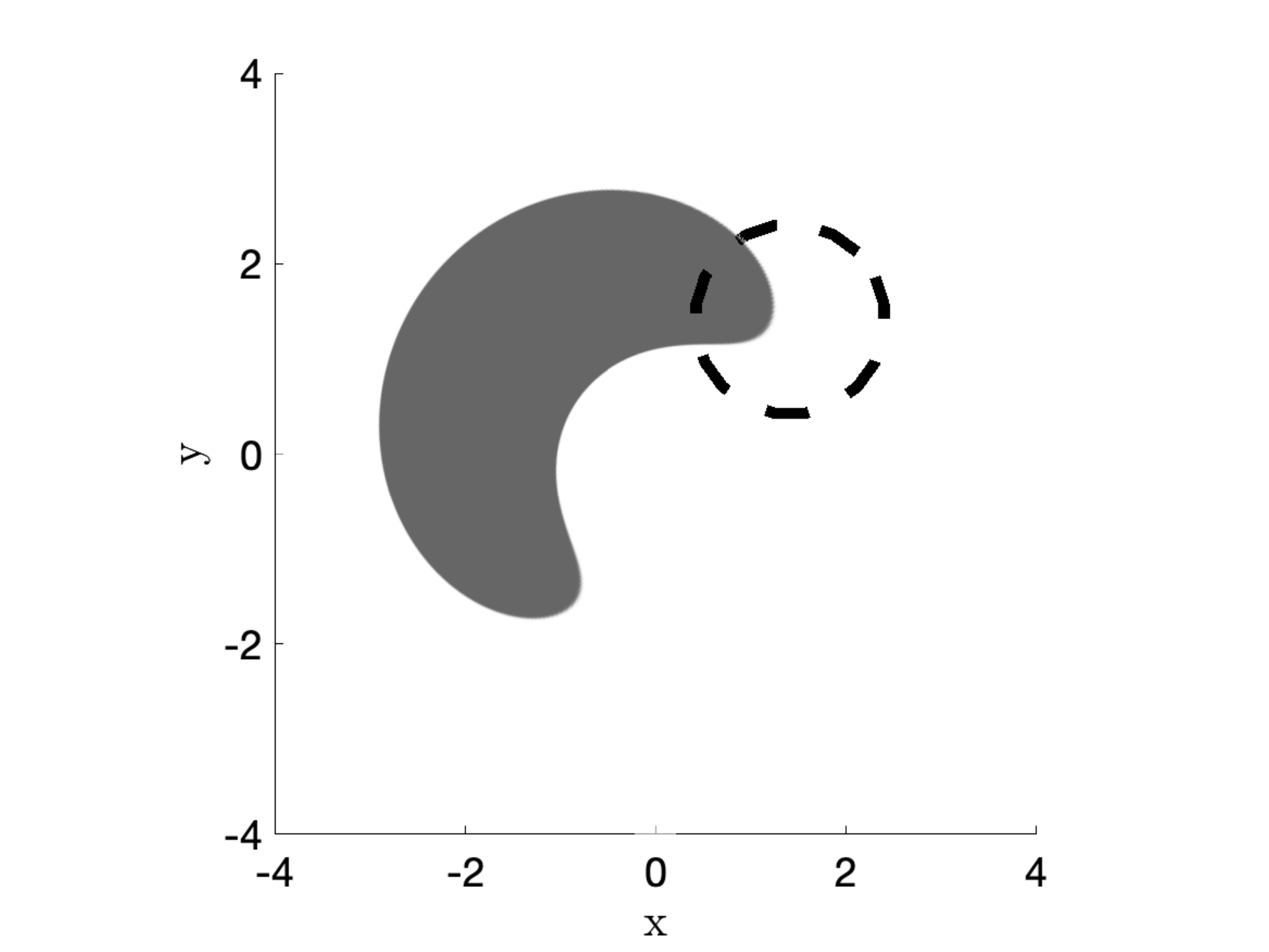}}
\end{center}
\end{minipage}
\end{minipage}
\begin{minipage}{\linewidth}
\vspace{8pt}
\begin{minipage}{0.49\linewidth}
\begin{center}
$t=1.25$ \\
\scalebox{0.32}{\includegraphics{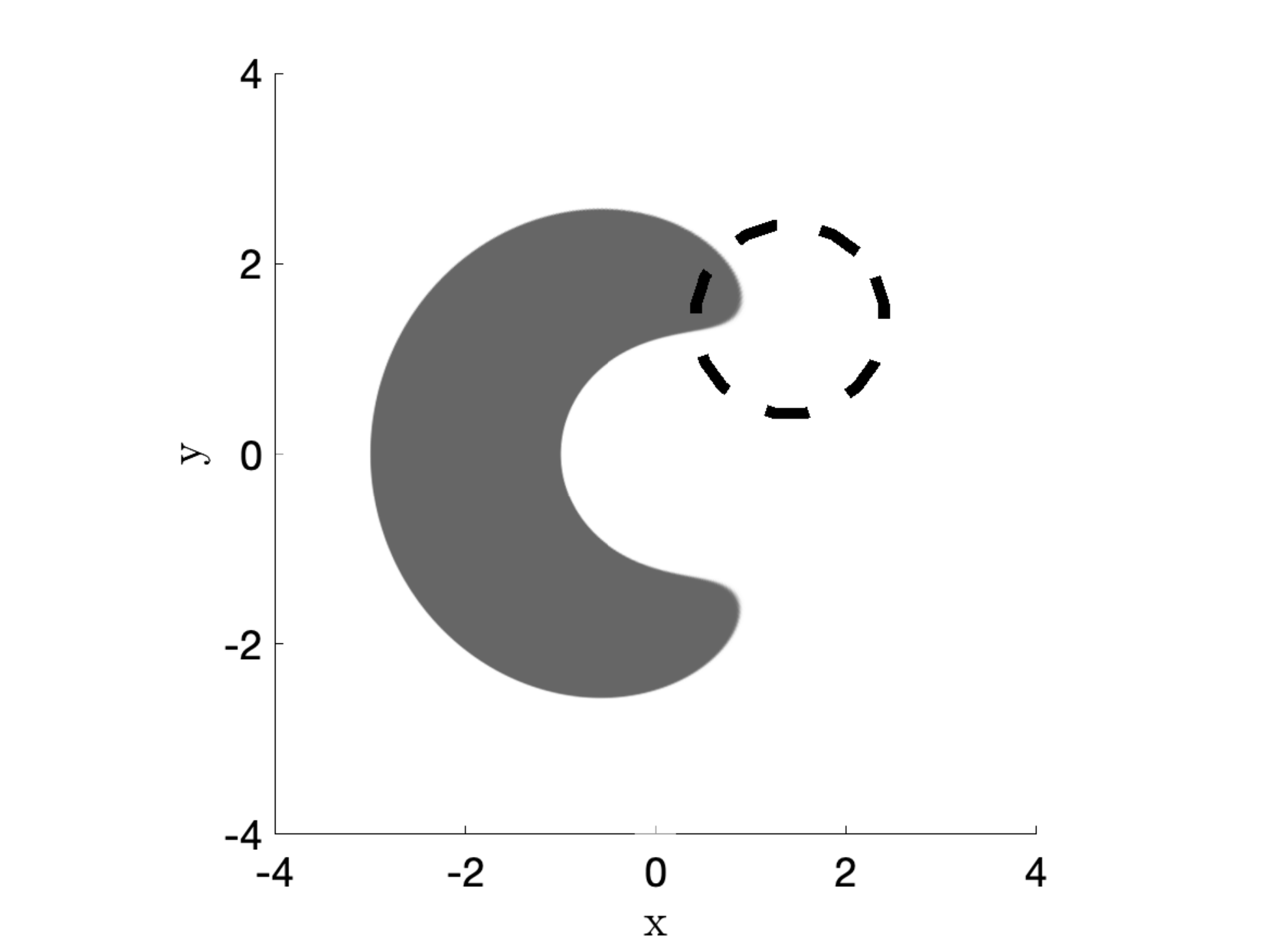}}
\end{center}
\end{minipage}
\begin{minipage}{0.49\linewidth}
\begin{center}
$t=2$ \\
\scalebox{0.32}{\includegraphics{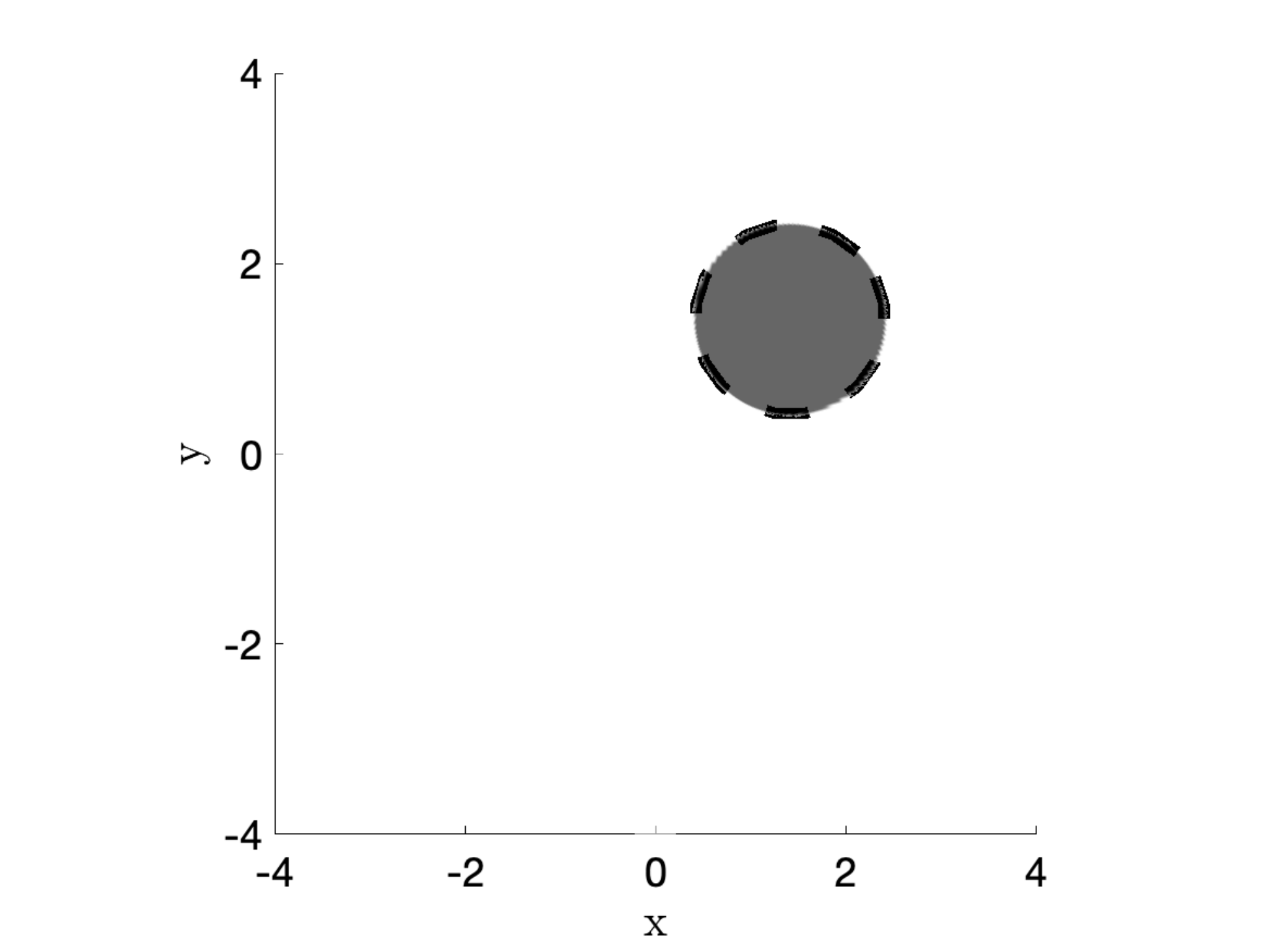}}
\end{center}
\end{minipage}
\end{minipage}
\begin{minipage}{\linewidth}
\vspace{8pt}
\begin{center}
\scalebox{0.39}{\includegraphics{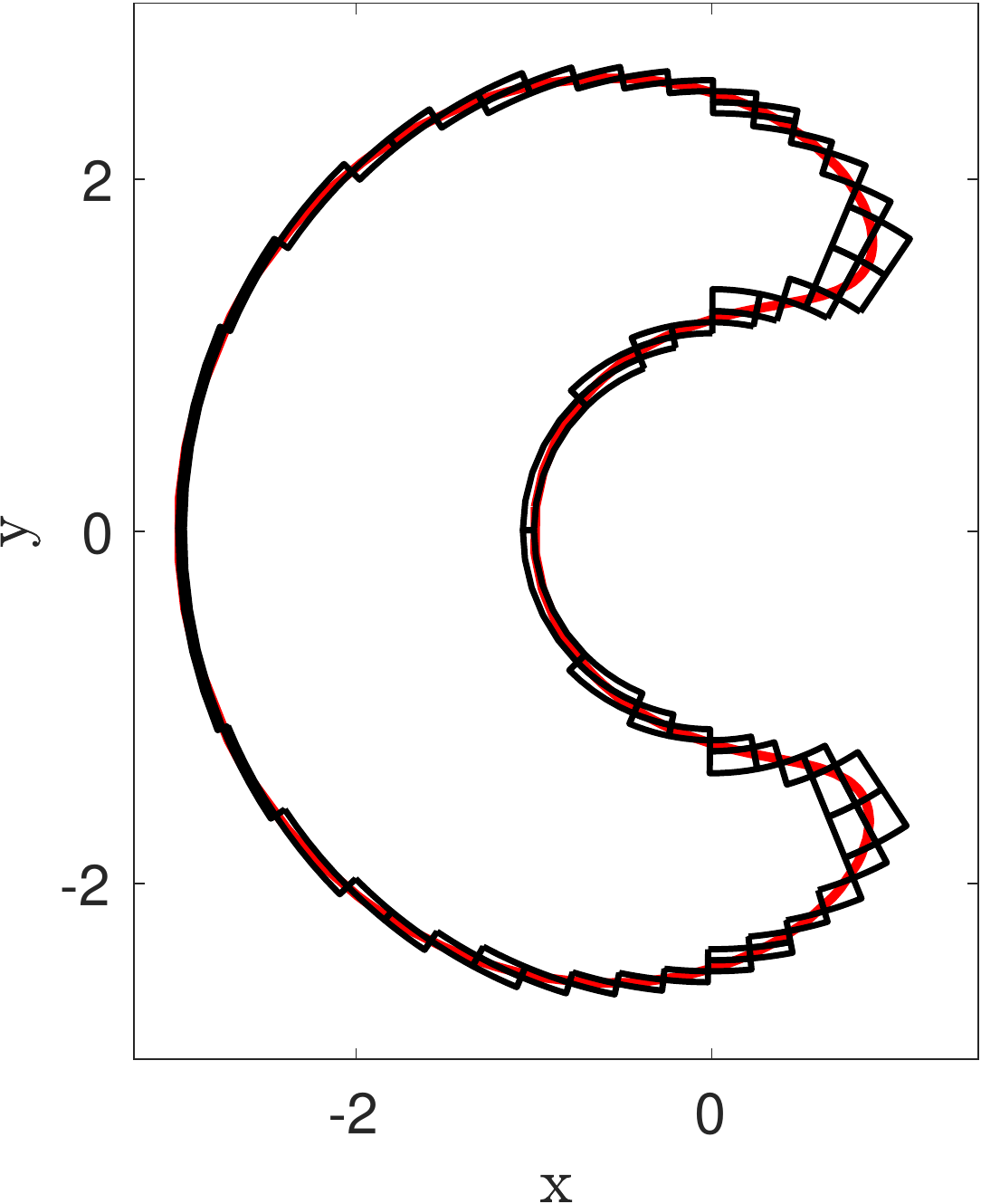}}
\hspace{50pt}
\scalebox{0.39}{\includegraphics{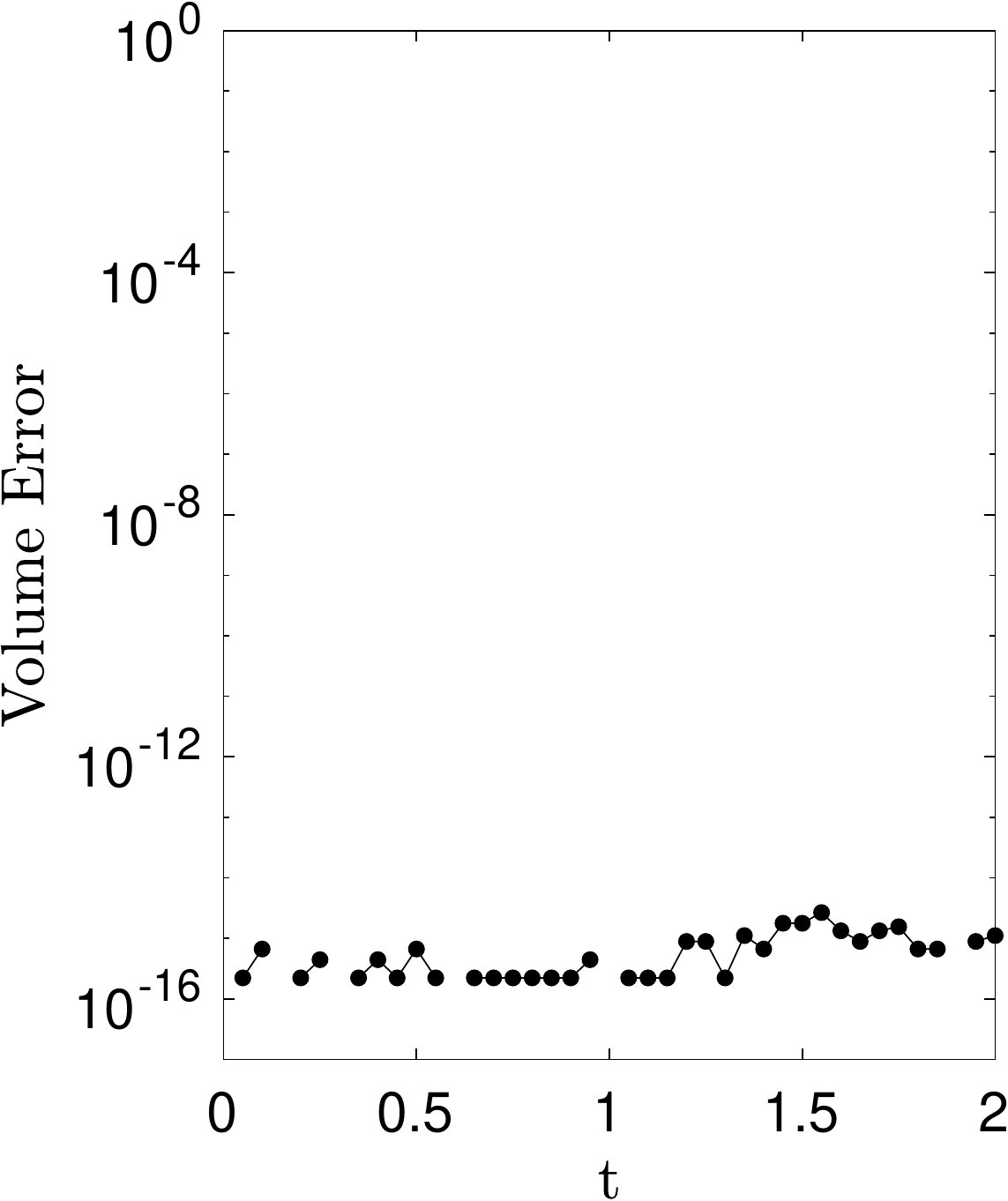}}
\end{center}
\end{minipage}
\caption{Stretching a unit disk under angular deformation. Dash line is the initial position at $t=0$. Top left is the position at $t=0.25$; Top right at $t=0.5$; Middle: $t=1.25$ (left) and finally at $t=2$ (right) back to original position. Bottom: a zoom-in view of the reconstructed interface (red) using patches with volume fractions between 7.5\% and 95\% at $t=1.25$. The volume error stays relatively flat at the order of machine precision throughout the simulation.}
\label{fig:numex6}
\end{figure}

Using the unit disk provided in Experiment 5, we can perform angular deformation of it in the $\theta$ direction only. We can similarly use the cosine multiplier as in Experiment 2 so that the disk is recovered at the original position after the deformation. \cref{fig:numex6} shows its dynamic in the cartesian coordinates recorded at different times. Finally, the disk is recovered at its original position, and the volume error remains at machine zero throughout. At $t=1.25$, we construct the approximation of the interface using the same technique done in Experiment 4. In addition to using only a centroid in a patch, to get more interface points nearby the patch, we can add an extra $1$ or $2$ more points along the direction of the tangent line, passing the centroids as initial guesses for the root-finding step. Having additional interface points, typically where the curvature changes quickly, will help us better parametrize the interface when needed. The bottom left figure of \cref{fig:numex6} shows the patches with volume fractions between 7.5\% and 95\%, and the red line is the reconstructed interface. 

\section{Discussion}
We should point out that the adaptive partition of unity method with moving patches described here is similar in flavor to the volume-of-fluid (VOF) method \cite{hirt_volume_1981}. We are not creating divergence-free bases (not explicitly on purpose) to approximate the two-medium function $f$. Instead, patches containing pieces of materials or volume fractions are moving along with the vector field; hence the volume conservation can be maintained with high accuracy. 

The method described here looks relatively simple for linear divergence-free types vector fields that affect material contour in a specific way. Moreover, rectangular patches can be made of different sizes in those cases, and the global interpolant can still be spectrally accurate. As provided in one of the numerical examples, one may switch to using polar domains for vector fields with pure rotations.

The patches may break one or more constraints provided in \cref{sec:movpatch} for problems that involve more complicated nonlinear vector fields. One possible remedy is to couple it with a method that domain decomposes the flow field regions so that locally, the linearized flow behaves like the ones provided in \cref{tab:typesvel}. In other words, the linearized model (Jacobian) regions of the nonlinear ODEs. We are investigating this approach for future study, especially comparing the proposed method with some standard benchmark problems found in \cite{sethian_level_1999, von_funck_vector_2006, ruuth2003simple, zhao_capturing_1998, stewart_improved_2008, takahashi_volume_2016, zolfaghari_simulations_2017, harvie_deformation_2008, dyadechko2005moment, tryggvason_direct_2011} for 2D and 3D cases. 

Two of our numerical experiments also showed that although the interface is not explicitly tracked, the highly accurate APU interpolant can be used for reconstructing the interface. Additionally, differentiation and integration operators can be done in the same framework. This is useful when one wants to utilize them to model material structure or fluid-structure interactions in multi-phase flow simulations. 
\section{Conclusions}
\label{sec:conclusions}

This study demonstrates that the adaptive partition of unity method with moving patches works very well with linear divergence-free vector fields while maintaining volume conservation with high accuracy. The method is simple to implement for shape deformation problems and in the same spirit as the VOF method but with spectrally accurate global interpolants. Using the technique for other types of vector fields or other more difficult benchmark problems is currently under investigation.

\section*{Acknowledgments}
We thank Kevin Aiton and Toby Driscoll for making the adaptive partition of unity codes open-source and UMass Dartmouth Center for Scientific Computing and Data Science Research (CSCDR) for providing rapid prototyping servers and computing facilities for this project. This work was funded by the National Science Foundations DMS-2012011. 

%Put references here.
\bibliographystyle{siamplain}
\bibliography{references}
\end{document}